\documentclass[final,onefignum,onetabnum]{siamonline171218}

\usepackage{lipsum}
\usepackage{amsfonts}
\usepackage{graphicx}
\usepackage{epstopdf}
\usepackage{algorithmic}
\usepackage{amssymb}
\usepackage{mathpazo}
\usepackage{mathrsfs}
\usepackage{amsmath}

\DeclareFontFamily{U}{mathx}{}
\DeclareFontShape{U}{mathx}{m}{n}{<-> mathx10}{}
\DeclareSymbolFont{mathx}{U}{mathx}{m}{n}
\DeclareMathAccent{\widehat}{0}{mathx}{"70}
\DeclareMathAccent{\widecheck}{0}{mathx}{"71}

\ifpdf
  \DeclareGraphicsExtensions{.eps,.pdf,.png,.jpg}
\else
  \DeclareGraphicsExtensions{.eps}
\fi
\newcommand{\dd}{\mathrm{d}}
\newcommand{\ee}{\mathrm{e}}
\newcommand{\ii}{\mathrm{i}}
\renewcommand{\varepsilon}{\epsilon}
\newcommand{\const}{\mathfrak{a}}
\def\R{{\mathbb R}}
\def\C{{\mathbb C}}

\def\build#1_#2^#3{\mathrel{
\mathop{\kern 0pt#1}\limits_{#2}^{#3}}}
\def\td_#1,#2{\mathrel{\mathop{\build\longrightarrow_{#1\rightarrow #2}^{}}}}

\def\e{\varepsilon}

\def\phe{\varphi}

\usepackage[shortlabels]{enumitem}
\setlist[enumerate]{leftmargin=.5in}
\setlist[itemize]{leftmargin=.5in}

\newsiamremark{remark}{Remark}
\newsiamremark{hypothesis}{Hypothesis}
\crefname{hypothesis}{Hypothesis}{Hypotheses}
\newsiamthm{claim}{Claim}
\newsiamthm{example}{Example}

\headers{The Benjamin-Ono Initial-Value Problem for Rational Data}{E. Blackstone, L. Gassot, P. G\'erard, and P. D. Miller}

\title{The Benjamin-Ono Initial-Value Problem for Rational Data with Application to Long-Time Asymptotics and Scattering
}

\author{Elliot Blackstone\thanks{Department of Mathematics, University of Michigan, Ann Arbor, MI 
  (\email{eblackst@umich.edu}, \url{http://www.math.lsa.umich.edu/}).}
\and Louise Gassot\thanks{CNRS and Department of Mathematics, University of Rennes, Rennes, France
	(\email{louise.gassot@cnrs.fr}).}
 \and
 Patrick G\'erard\thanks{Laboratoire de Mathématiques d'Orsay,  Université Paris-Saclay, Orsay, France (\email{patrick.gerard@universite-paris-saclay.fr}).}
\and Peter D. Miller\thanks{Department of Mathematics, University of Michigan, Ann Arbor, MI 
  (\email{millerpd@umich.edu}).}
}

\usepackage{amsopn}

\makeatletter
\newcommand*{\addFileDependency}[1]{
  \typeout{(#1)}
  \@addtofilelist{#1}
  \IfFileExists{#1}{}{\typeout{No file #1.}}
}
\makeatother

\ifpdf
\hypersetup{
  pdftitle={The Benjamin-Ono Initial-Value Problem for Rational Data},
  pdfauthor={Elliot Blackstone and Louise Gassot and Patrick Gérard and Peter D. Miller}
}
\fi

\usepackage{comment}
\usepackage{tikz}
\usetikzlibrary{arrows}
\usetikzlibrary{decorations.pathmorphing}
\usetikzlibrary{decorations.markings}
\usetikzlibrary{patterns}
\usetikzlibrary{automata}
\usetikzlibrary{positioning}
\usepackage{tikz-cd}

\tikzset{->-/.style={decoration={markings,mark=at position #1 with {\arrow[thick]{>}}},postaction={decorate}}}
	
\tikzset{-<-/.style={decoration={markings,mark=at position #1 with{\arrow[thick]{<}}},postaction={decorate}}}

\tikzset{
  frac arrow/.style={postaction={decorate,decoration={
        markings,
        mark=at position #1 with {\arrow[]{Stealth[round]}}
      }}},
}

\begin{document}

\maketitle

\begin{abstract}
  We show that the initial-value problem for the Benjamin-Ono equation on $\mathbb{R}$ with $L^2(\mathbb{R})$ rational initial data with only simple poles can be solved in closed form via a determinant formula involving contour integrals.  The dimension of the determinant depends on the number of simple poles of the rational initial data only and the matrix elements depend explicitly on the independent variables $(t,x)$ and the dispersion coefficient $\epsilon$.  This allows for various interesting asymptotic limits to be resolved quite efficiently. As an example, and as a first step towards establishing the soliton resolution conjecture, we prove that the solution with initial datum equal to minus a soliton exhibits scattering.
\end{abstract}

\begin{keywords}
  Benjamin-Ono equation, Cauchy problem, exact solutions
\end{keywords}

\begin{AMS}
35C05, 35Q51, 37K10
\end{AMS}

\section{Introduction}
The Benjamin-Ono equation
\begin{equation}
\partial_t u + \partial_x(u^2) = \epsilon \partial_x|D_x|u,\quad x\in\mathbb{R},\quad t\in\mathbb{R},\quad \epsilon>0,
\label{eq:BO}
\end{equation}
where $\epsilon>0$ measures the strength of the dispersion, is a well-known asymptotic model for internal waves in stratified fluids that was first proposed more than 50 years ago \cite{Benjamin66,Benjamin67}.  The dispersive term is either defined by the Fourier multiplier $\widehat{|D_x|u}(\xi)=|\xi|\widehat{u}(\xi)$  or by the singular integral
\begin{equation}|D_x|f(x):=-\frac{1}{\pi}{\mathrm{P.V.}}\int_\mathbb{R}\frac{f'(y)\,\dd y}{y-x}.
\end{equation}
Solutions $u=u(t,x)$ represent the vertical displacement of an interface (the \emph{pycnocline}) between fluid layers of differing density, with the lower, denser layer being assumed to be infinitely deep.  

It has been known since the late 1970's and early 1980's that the Benjamin-Ono equation has many of the features of a completely integrable dispersive nonlinear wave equation, the prototype of which is the Korteweg-de Vries equation.  Indeed, the Benjamin-Ono equation has a B\"acklund transformation and infinitely-many conservation laws \cite{Nakamura79}, can be represented as the compatibility condition of a Lax pair \cite{BockKruskal79}, is part of an infinite hierarchy of commuting Hamiltonian flows \cite{FokasFuchssteiner81}, has a variety of exact solutions including multi-solitons \cite{Matsuno79} and multi-phase waves \cite{DobrokhotovK91}, and has been solved by a formal inverse-scattering transform method \cite{FokasAblowitz83,KaupMatsuno98}.  A generalization of the latter approach, the method of commuting flows, has recently been used to prove a sharp well-posedness \cite{KillipLaurensVisan23-BO} result for the Benjamin-Ono equation.

While the inverse-scattering transform method in its original form has yet to be made fully rigorous (see \cite{KleinSaut21,Wu16,Wu17}), there have been some developments in its use to study the initial-value problem on the line $x\in\mathbb{R}$ with initial data $u_0$ in a suitable space of functions with decay at infinity.  Generalizing from an illustrative example in the paper \cite{KodamaAblowitzSatsuma82}, it was shown in \cite{MillerW16a} that all of the scattering data in the Fokas-Ablowitz inverse-scattering transform can be obtained essentially in explicit form when $u_0$ is a rational function.  For instance, the discrete eigenvalues of the Lax operator are the zeros of an Evans function that is given explicitly as a determinant, the elements of which are analytic functions defined by contour integrals involving $u_0$, and the dimension of the determinant is precisely the number of poles of $u_0$ in the upper complex half-plane.  The reflection coefficient also has an explicit expression in terms of the solution of a linear system of the same dimension with coefficients that are similar contour integrals.  Another remarkable recent result \cite{Gerard22} is that the formal inverse-scattering algorithm can be ``closed'' to yield a formula for the solution $u(t,x)$ of the initial-value problem with completely general initial data $u_0\in L^2(\mathbb{R})\cap L^{\infty}(\mathbb{R})$ in terms of the resolvent of a linear operator depending explicitly on $u_0$ (see \eqref{eq:Gerard} below). This result has been extended to initial data $u_0 \in L^2(\mathbb{R})$ in~\cite{Chen24} by using dispersive properties of the free Schr\"odinger evolution.  Although the steps of the inverse-scattering transform themselves have not been fully justified, the resulting closed formula for $u(t,x)$ has been proved to provide the solution of the initial-value problem nonetheless.  

The aim of this paper is to apply the exact solution formula in \cite{Gerard22} to the special case that $u_0\in L^2(\mathbb{R})$ is a real-valued rational function with simple poles:
\begin{equation}
u(0,x)=u_0(x):=\sum_{n=1}^N\left(\frac{c_n}{x-p_n}+\frac{{c_n^*}}{x-{p_n^*}}\right),
\quad \operatorname{Im}(p_n)>0.
\label{eq:rationalIC}
\end{equation}
Clearly $u_0\in L^2(\mathbb{R})$ for all pole locations $\{p_n\}$ in the upper half-plane and for all complex values of the coefficients $\{c_n\}$.  We will show how completely explicit the formula from \cite{Gerard22} becomes for such rational data.  In particular, the solution $u(t,x)$ of \eqref{eq:BO} for initial data $u(0,x)=u_0(x)$ given by \eqref{eq:rationalIC} turns out to be given by a ratio of $(N+1)\times (N+1)$ determinants, the entries of which are contour integrals involving the initial data.  This result is therefore not unlike the computation in \cite{MillerW16a} of the Fokas-Ablowitz scattering data for $u_0$ of the form \eqref{eq:rationalIC}, but more to the point it also completely solves the corresponding inverse problem and returns explicitly the solution $u(t,x)$ instead of just its data in the transform domain.  Obviously the availability of such a formula (see Theorem~\ref{thm:inversion-formula} below) makes it immediately tractable to analyze the solution $u(t,x)$ in various interesting asymptotic limits. We have successfully used it to study in particular the small-dispersion limit $\epsilon\to 0$ for fixed rational initial data $u_0$ to exhibit multi-phase wave asymptotics away from wave caustic curves in the $(t,x)$-plane \cite{SmallDispersionBO-1} as well as to prove universality theorems \cite{SmallDispersionBO-2} in the vicinity of such curves and their intersection points. It is also very useful  to analyze the long-time limit \cite{LongTime-BO} and we believe it will allow us to prove the soliton resolution conjecture, in which we expect that the radiation term exhibits scattering if $\langle x\rangle u_0\in L^2$, and modified scattering otherwise. As an example, we use the new exact solution formula to prove that this is the case for an initial datum equal to minus a soliton; see Theorem~\ref{thm:minussoliton} below.

Now we explain our main results.  
Choosing $\log$ to be the principal branch such that $|\operatorname{Im}(\log(\cdot))|<\pi$, we first define for $y\in\mathbb{R}$ negative and sufficiently large in absolute value,
\begin{equation}
h(y):=\frac{1}{4t}(y-x)^2+\sum_{n=1}^N\left[c_n\log(
y-p_n)+c_n^*\log(y-p_n^*)\right],
\label{eq:h-def}
\end{equation}
so that comparing with \eqref{eq:rationalIC},
\begin{equation}
    h'(y)=\frac{y-x}{2t}+u_0(y).
 \label{eq:hprime}
\end{equation}

Note that when $y<0$ with $|y|$ sufficiently large, $h(y)\in\mathbb{R}$.  
We want to analytically continue $y\mapsto h(y)$ to a maximal domain that is generally more complicated than implied by using the principal branch of the logarithm in \eqref{eq:h-def} and taking $y=z$ to be complex.  For this purpose, we start by allowing for general branch cuts $\{\Gamma_n,\bar{\Gamma}_n\}_{n=1}^N$ assumed only to have the following properties.
\begin{definition}[Branch cuts of $h$]
\label{def:B}
The branch cuts $\Gamma_1,\dots,\Gamma_N$ are pairwise disjoint piecewise-smooth curves each emanating from exactly one of the poles $\{p_n\}_{n=1}^N$ and tending to $z=\infty$ in the direction asymptotic to the ray $\arg(z)=3\pi/4$. All of these branch cuts are assumed to lie in a half-plane $\mathrm{Im}(z)>-\delta$ for some $\delta>0$ sufficiently small (in particular, we assume $\delta<\min_n\{\mathrm{Im}(p_n)\}$.

The branch cuts $\bar{\Gamma}_1,\dots,\bar{\Gamma}_N$ are straight horizontal rays each emanating from exactly one of the conjugate poles $\{p_n^*\}_{n=1}^N$ and extending to $z=\infty$ in the left half-plane.

\end{definition}
Given the branch cuts, we index the poles $\{p_n\}_{n=1}^N$ such that in the vicinity of $z=\infty$ in the upper half-plane, the branch cuts $\Gamma_1,\dots,\Gamma_N$ are ordered left-to-right.  See Figure~\ref{fig:BO-BasicContours}, left-hand panel.
We hence obtain a well-defined function $z\mapsto h(z)$ by analytic continuation from large negative real values of $z$ where $h(z)$ is given by \eqref{eq:h-def} to the domain \begin{equation}
z\in\mathbb{C}\setminus (\Gamma_1\cup\cdots\cup \Gamma_N\cup \bar{\Gamma}_1\cup\cdots\cup\bar{\Gamma}_N).
\end{equation}
We next define some relevant integration contours in the $z$-plane, see Figure~\ref{fig:BO-BasicContours}, right-hand panel.

\begin{definition}[Contours]\label{def:C}
  Let $C_n$, $n=1,\dots,N$ denote the contour defined by one of the following alternatives:
\begin{itemize}
    \item If $\ii c_n/\epsilon$ is a strictly {negative} integer, then $C_n$ originates at $z=\infty$ in the direction $\arg(z)={3\pi/4}$ to the {left} of all branch cuts of $h$, lies in the domain of analyticity of $h$, and terminates at $z=p_n$.
    We call such an index $n$ \emph{exceptional}.
    \item Otherwise,  $C_n$ originates and terminates at $z=\infty$ in the direction $\arg(z)={3\pi/4}$, lies in the domain of analyticity of $h(z)$, and encircles with counterclockwise orientation precisely the branch cuts of $h(z)$ emanating from each of the points $z=p_m$, $1\le m\le n$.  Such an index $n$ will be called \emph{non-exceptional}.
\end{itemize}
Finally, we let $C_0$ denote a path in the domain of analyticity of $h(z)$ originating at $z=\infty$ in the direction $\arg(z)={3\pi/4}$  to the {left} of all  branch cuts $\Gamma_1,\dots,\Gamma_N$ of $h(z)$ and terminating at $z=\infty$ in the direction $\arg(z)=-{\pi/4}$.
\end{definition}

\begin{figure}
\hfill\begin{tikzpicture}
\node (p1) at (-1.2cm,0.5cm) [thick,black,circle,fill=lightgray,draw=black,inner sep=1.2pt]{};
\node (p2) at (1.25cm,0.65cm) [thick,black,circle,fill=lightgray,draw=black,inner sep=1.2pt]{}; 
\node (p3) at (0.3cm,0.9cm) [thick,black,circle,fill=lightgray,draw=black,inner sep=1.2pt]{}; 
\node (p4) at (-0.7cm,1.3cm) [thick,black,circle,fill=lightgray,draw=black,inner sep=1.2pt]{}; 
\node (p5) at (1.5cm,1.1cm) [thick,black,circle,fill=lightgray,draw=black,inner sep=1.2pt]{}; 
\node at (-0.9cm,0.5cm) {$p_1$};
\node at (1.55cm,0.65cm) {$p_2$};
\node at (0.6cm,0.9cm) {$p_3$};
\node at (-1.0cm,1.4cm) {$p_4$};
\node at (1.8cm,1.1cm) {$p_5$};
\node (p1s) at (-1.2cm,-0.5cm) [thick,black,circle,fill=lightgray,draw=black,inner sep=1.2pt]{};
\node (p2s) at (1.25cm,-0.65cm) [thick,black,circle,fill=lightgray,draw=black,inner sep=1.2pt]{}; 
\node (p3s) at (0.3cm,-0.9cm) [thick,black,circle,fill=lightgray,draw=black,inner sep=1.2pt]{}; 
\node (p4s) at (-0.7cm,-1.3cm) [thick,black,circle,fill=lightgray,draw=black,inner sep=1.2pt]{}; 
\node (p5s) at (1.5cm,-1.1cm) [thick,black,circle,fill=lightgray,draw=black,inner sep=1.2pt]{}; 
\node at (-0.9cm,-0.4cm) {$p_1^*$};
\node at (1.55cm,-0.65cm) {$p_2^*$};
\node at (0.6cm,-0.875cm) {$p_3^*$};
\node at (-0.4cm,-1.4cm) {$p_4^*$};
\node at (1.8cm,-1.1cm) {$p_5^*$};
\draw [->,>={Stealth[round]},thick,lightgray] (-1.6cm,0cm)--(2.2cm,0cm) {};
\draw [->,>={Stealth[round]},thick,lightgray] (0cm,-2cm) -- (0cm,2cm) {};
\draw [thick,red,decorate,
decoration={snake,amplitude=.2mm,segment length=1mm,post length=0mm}] (p1) -- (-1.6cm,0.9cm);
\draw [thick,red,decorate,
decoration={snake,amplitude=.2mm,segment length=1mm,post length=0mm}] (p5) -- (1cm,1.6cm);
\draw [thick,red,decorate,
decoration={snake,amplitude=.2mm,segment length=1mm,post length=0mm}] (p2) -- (0.5cm,-0.1cm) -- (0.0cm,-0.1cm) -- (-1.6cm,1.5cm);
\draw [thick,red,decorate,
decoration={snake,amplitude=.2mm,segment length=1mm,post length=0mm}] (p4) -- (0.4cm,0.2cm) -- (1.1cm,0.9cm) -- (0.4cm,1.6cm);
\draw [thick,red,decorate,
decoration={snake,amplitude=.2mm,segment length=1mm,post length=0mm}] (p3) -- (-0.4cm,1.6cm);
\draw [thick,red,decorate,
decoration={snake,amplitude=.2mm,segment length=1mm,post length=0mm}] (p1s) -- (-1.6cm,-0.5cm);
\draw [thick,red,decorate,
decoration={snake,amplitude=.2mm,segment length=1mm,post length=0mm}] (p2s) -- (-1.6cm,-0.65cm);
\draw [thick,red,decorate,
decoration={snake,amplitude=.2mm,segment length=1mm,post length=0mm}] (p3s) -- (-1.6cm,-0.9cm);
\draw [thick,red,decorate,
decoration={snake,amplitude=.2mm,segment length=1mm,post length=0mm}] (p4s) -- (-1.6cm,-1.3cm);
\draw [thick,red,decorate,
decoration={snake,amplitude=.2mm,segment length=1mm,post length=0mm}] (p5s) -- (-1.6cm,-1.1cm);

\node at (-1.8cm,0.9cm) {$\color{red}\Gamma_1$};
\node at (-1.8cm,1.6cm) {$\color{red}\Gamma_2$};
\node at (-0.5cm,1.8cm) {$\color{red}\Gamma_3$};
\node at (0.3cm,1.8cm) {$\color{red}\Gamma_4$};
\node at (0.9cm,1.8cm) {$\color{red}\Gamma_5$};

\node at (-1.8cm,-0.5cm) {$\color{red}\bar{\Gamma}_1$};
\node at (-2.2cm,-0.7cm) {$\color{red}\bar{\Gamma}_2$};
\node at (-1.8cm,-0.9cm) {$\color{red}\bar{\Gamma}_3$};
\node at (-2.2cm,-1.1cm) {$\color{red}\bar{\Gamma}_5$};
\node at (-1.8cm,-1.3cm) {$\color{red}\bar{\Gamma}_4$};

\draw [thick,black,dashed] (-1.6cm,-0.425cm) -- (2.2cm,-0.425cm);
\draw [thick,black,<->,>={Stealth[round]}] (1.9cm,0cm) -- (1.9cm,-0.425cm);
\node at (2.1cm,-0.2cm) {$\delta$};
\end{tikzpicture} \hfill%
\begin{tikzpicture}
\node (p1) at (-1.2cm,0.5cm) [thick,black,circle,fill=lightgray,draw=black,inner sep=1.2pt]{};
\node (p2) at (1.25cm,0.65cm) [thick,black,circle,fill=lightgray,draw=black,inner sep=1.2pt]{}; 
\node (p3) at (0.3cm,0.9cm) [thick,black,circle,fill=lightgray,draw=black,inner sep=1.2pt]{}; 
\node (p4) at (-0.7cm,1.3cm) [thick,black,circle,fill=lightgray,draw=black,inner sep=1.2pt]{}; 
\node (p5) at (1.5cm,1.1cm) [thick,black,circle,fill=lightgray,draw=black,inner sep=1.2pt]{}; 
\node (p1s) at (-1.2cm,-0.5cm) [thick,black,circle,fill=lightgray,draw=black,inner sep=1.2pt]{};
\node (p2s) at (1.25cm,-0.65cm) [thick,black,circle,fill=lightgray,draw=black,inner sep=1.2pt]{}; 
\node (p3s) at (0.3cm,-0.9cm) [thick,black,circle,fill=lightgray,draw=black,inner sep=1.2pt]{}; 
\node (p4s) at (-0.7cm,-1.3cm) [thick,black,circle,fill=lightgray,draw=black,inner sep=1.2pt]{}; 
\node (p5s) at (1.5cm,-1.1cm) [thick,black,circle,fill=lightgray,draw=black,inner sep=1.2pt]{}; 
\draw [->,>={Stealth[round]},thick,lightgray] (-1.6cm,0cm)--(2.2cm,0cm) {};
\draw [->,>={Stealth[round]},thick,lightgray] (0cm,-2cm) -- (0cm,2cm) {};
\draw [thick,red,decorate,
decoration={snake,amplitude=.2mm,segment length=1mm,post length=0mm}] (p1) -- (-1.6cm,0.9cm);
\draw [thick,red,decorate,
decoration={snake,amplitude=.2mm,segment length=1mm,post length=0mm}] (p5) -- (1cm,1.6cm);
\draw [thick,red,decorate,
decoration={snake,amplitude=.2mm,segment length=1mm,post length=0mm}] (p2) -- (0.5cm,-0.1cm) -- (0.0cm,-0.1cm) -- (-1.6cm,1.5cm);
\draw [thick,red,decorate,
decoration={snake,amplitude=.2mm,segment length=1mm,post length=0mm}] (p4) -- (0.4cm,0.2cm) -- (1.1cm,0.9cm) -- (0.4cm,1.6cm);
\draw [thick,red,decorate,
decoration={snake,amplitude=.2mm,segment length=1mm,post length=0mm}] (p3) -- (-0.4cm,1.6cm);
\draw [thick,red,decorate,
decoration={snake,amplitude=.2mm,segment length=1mm,post length=0mm}] (p1s) -- (-1.6cm,-0.5cm);
\draw [thick,red,decorate,
decoration={snake,amplitude=.2mm,segment length=1mm,post length=0mm}] (p2s) -- (-1.6cm,-0.65cm);
\draw [thick,red,decorate,
decoration={snake,amplitude=.2mm,segment length=1mm,post length=0mm}] (p3s) -- (-1.6cm,-0.9cm);
\draw [thick,red,decorate,
decoration={snake,amplitude=.2mm,segment length=1mm,post length=0mm}] (p4s) -- (-1.6cm,-1.3cm);
\draw [thick,red,decorate,
decoration={snake,amplitude=.2mm,segment length=1mm,post length=0mm}] (p5s) -- (-1.6cm,-1.1cm);

\draw [thick,postaction={frac arrow=0.85}] (-1.6cm,0.7cm) -- (-1.2cm,0.3cm) -- (-1.0cm,0.5cm) -- (-1.6cm,1.1cm); 
\draw [thick,postaction={frac arrow=0.925}] (-1.6cm,0.65cm) -- (-0.75cm,-0.2cm) -- (1.25cm,-0.2cm) -- (p2); 
\draw [thick,postaction={frac arrow=0.97}] (-1.6cm,0.6cm) -- (-0.75cm,-0.25cm) -- 
(1.6cm,-0.25cm) -- (1.6cm,0.6cm) -- (1.25cm,0.95cm) -- (0.35cm,0.05cm) -- (-0.9cm,1.3cm) -- (-0.7cm,1.5cm) -- (0.2cm,0.6cm) -- (0.5cm,0.9cm) -- (-0.2cm,1.6cm); 
\draw [thick,postaction={frac arrow=0.95}] (-1.6cm,0.55cm) -- (-0.75cm,-0.3cm) -- (1.7cm,-0.3cm) -- (1.7cm,0.6cm) -- (0.7cm,1.6cm); 
\draw [thick,postaction={frac arrow=0.95}] (-1.6cm,0.5cm) -- (-0.75cm,-0.35cm) -- (1.8cm,-0.35cm) -- (1.8cm,1.1cm) -- (1.3cm,1.6cm); 
\draw [thick,postaction={frac arrow=0.95}] (-1.6cm,0.45cm) -- (-0.75cm,-0.4cm) -- (1.8cm,-0.4cm) -- (2.2cm,-0.8cm); 
\node at (-1.8cm,1.1cm) {$C_1$};
\node at (1.0cm,0.025cm) {$C_2$};
\node at (-0.2cm,1.8cm) {$C_3$};
\node at (0.7cm,1.8cm) {$C_4$};
\node at (1.3cm,1.8cm) {$C_5$};
\node at (1.8cm,-0.7cm) {$C_0$};
\draw [thick,black,dashed] (-1.6cm,-0.425cm) -- (2.2cm,-0.425cm);
\end{tikzpicture}\hfill%
    \caption{Left:  admissible branch cuts of $h(z)$ in the $z$-plane for a rational initial condition with $N=5$.  Right:  corresponding contours $C_1,C_2,C_3,C_4,C_5$ and $C_0$ for a situation where $2$ is the only exceptional index.}
\label{fig:BO-BasicContours}
\end{figure}

Let $\mathbf{A}(t,x),\mathbf{B}(t,x)$ be two $(N+1)\times (N+1)$ matrices defined for $\epsilon>0$ and $t>0$ by 
\begin{align}
A_{j1}(t,x):=\int_{C_{j-1}}u_0(z)\ee^{-\ii h(z)/\epsilon}\,\dd z, 
	&\quad
A_{jk}(t,x):=\int_{C_{j-1}}\frac{\ee^{-\ii h(z)/\epsilon}\,\dd z}{z-p_{k-1}},
\label{eq:Amatrix}
\\
B_{j1}(t,x):=\int_{C_{j-1}}\ee^{-\ii h(z)/\epsilon}\,\dd z,
	&\quad
B_{jk}(t,x):=\int_{C_{j-1}}\frac{\ee^{-\ii h(z)/\epsilon}\,\dd z}{z-p_{k-1}}=A_{jk}(t,x),
\label{eq:Bmatrix}
\end{align}
where the indices satisfy $1\leq j\leq N+1$ and $ 2\leq k\leq N+1$.
Note that if $j$ is an exceptional index, then $\ee^{-\ii h(z)/\epsilon}$ is actually a single-valued analytic function of $z$ on a neighborhood of the branch cut $\Gamma_j$ of $h(z)$, and this function vanishes to at least first order at the branch point $p_j$.  Finally, we introduce a related matrix $\overline{\mathbf{B}}(t,x)$, the first column of which is the same as that of $\mathbf{B}(t,x)$ while 
\begin{equation}\label{eq:Btilde}
    \overline{B}_{jk}=\ee^{\ii(x-p_{k-1})^2/(4t\epsilon)}B_{jk}=\ee^{\ii(x-p_{k-1})^2/(4t\epsilon)}\int_{C_{j-1}}\frac{\ee^{-\ii h(z)/\epsilon}}{z-p_{k-1}}\,\dd z,\quad k=2,\dots,N+1.
\end{equation}

\begin{theorem}[Solution of Benjamin-Ono for rational initial data]\label{thm:inversion-formula}
Let $\epsilon>0$.  The solution of the Cauchy initial-value problem for the Benjamin-Ono equation \eqref{eq:BO} with rational initial condition $u(x,0)=u_0(x)$ of the form \eqref{eq:rationalIC} is 
\begin{equation}
u(t,x)=\Pi u(t,x) +
\Pi u(t,x)^*=2\mathrm{Re}(\Pi u(t,x)),\quad t>0,
\label{eq:u-Piu}
\end{equation}
\begin{equation}
\Pi u(t,x)
	=\frac{\det(\mathbf{A}(t,x))}{\det(\mathbf{B}(t,x))}.
\label{eq:lambda-formula}
\end{equation}
An equivalent formula is
\begin{equation}
   \Pi  u(t,x)=\ii \epsilon\frac{\partial}{\partial x}\log(\det(\overline{\mathbf{B}}(t,x))).
\label{eq:tau-form}
\end{equation}
Also, we have $\det(\mathbf{B}(t,x))\neq 0$ and $\det(\overline{\mathbf{B}}(t,x))\neq 0$ for all $(t,x)\in\mathbb{R}^2$ with $t>0$ and all $\epsilon>0$.
\end{theorem}
The expression $\Pi u$ denotes the action of a Cauchy-Szeg\H{o} projection with respect to the argument $x$, which is explained in Section~\ref{section:inversion}.

\begin{remark}[The limit $t\downarrow 0$ and negative times]
The elements of the matrices $\mathbf{A}(t,x)$ and $\mathbf{B}(t,x)$, and hence also $\Pi u(t,x)$ given in \eqref{eq:lambda-formula}, are undefined for $t=0$ due to the factor of $t^{-1}$ in the expression \eqref{eq:h-def} for $h(z)$.  However, since $h(z)$ appears in the exponent of each integrand, a steepest-descent analysis shows that as $t\downarrow 0$, we have $\Pi u(t,x)\to \Pi u_0(x)$ for each $x\in\mathbb{R}$.  In this way, the solution can be defined by continuity for $t=0$.  Furthermore, to continue the solution to $t<0$ we may use the fact that 
if $u$ is a solution to~\eqref{eq:BO}, then $(t,x)\mapsto u(-t,-x)$ is also a solution to~\eqref{eq:BO}.  Note that the map $u_0(x)\mapsto u_0(-x)$ preserves the class of rational functions in $L^2(\mathbb{R})$ and amounts to the change of rational data $\{p_n\}_{n=1}^N\mapsto \{-p_n^*\}_{n=1}^N$ and $\{c_n\}_{n=1}^N\mapsto\{-c_n^*\}_{n=1}^N$ in \eqref{eq:rationalIC}.
\label{rem:tzero-and-negative}
\end{remark}

\begin{remark}[Choice of branch cuts]
    Since it is known that the Cauchy problem with data $u_0\in H^2(\mathbb{R})$ has a unique solution~\cite{Saut79}, the real part of the ratio of determinants in~\eqref{eq:lambda-formula} is insensitive to the choice of admissible branch cuts and integration contours consistent with Definitions~\ref{def:B} and~\ref{def:C}.  For some applications it is sufficient to choose straight-ray cuts for which $\Gamma_n$ is simply $p_n+\ee^{3\pi\ii/4}\mathbb{R}_+$.  Furthermore, the branch cuts $\bar{\Gamma}_j$ in the lower half-plane can be chosen to be quite arbitrary as long as they do not intersect any of the contours of integration $C_n$.  The flexibility of choosing general contours is especially useful in asymptotic problems such as those considered in \cite{SmallDispersionBO-1} and \cite{SmallDispersionBO-2}.
\label{rem:straight-rays}
\end{remark}

\begin{remark}[Two equivalent formul\ae]
    Of the two equivalent formul\ae\ for $u(t,x)$, \eqref{eq:lambda-formula} is preferable in situations where one is interested in the analysis of $u(t,x)$ in some asymptotic limit such as $\epsilon\to 0$ or $t\to\infty$, because one does not need to determine the effect of differentiation on error terms arising by expanding the determinants of $\mathbf{A}(t,x)$ and $\mathbf{B}(t,x)$.  On the other hand, the formula \eqref{eq:tau-form} is more compact, better for exact computations such as those in Appendix~\ref{section:soliton} below, 
    and also is reminiscent of formul\ae\ for more traditional exact solutions as described in Remark~\ref{rem:remarkable}.
\end{remark}

 \begin{remark}[Comparison with the Calogero-Moser equation]
Regarding the Calogero-Moser derivative NLS equation studied in ~\cite{GerardLenzmann}, a similar formula as~\eqref{eq:lambda-formula} for rational initial data can be derived from the general solution formula  established in~\cite{KillipLaurensVisan23}. However, no formulation comparable to~\eqref{eq:tau-form} is known. This could be explained by the different nature between the two equations. In fact, the solutions of the Benjamin-Ono equation stay bounded with respect to every Sobolev norm~\cite{Nakamura79,BockKruskal79}. On the contrary, for the Calogero-Moser derivative NLS equation, solutions with initial data of mass greater than the mass of the ground state can exhibit turbulent behavior~\cite{HoganKowalski24} or even blow-up~\cite{KimKimKwon24}.
\end{remark}

\begin{remark}[Comparison with the $N$-soliton formula]
    The exact solution formula asserted in Theorem~\ref{thm:inversion-formula} for the Benjamin-Ono initial-value problem with rational initial data $u_0(x)$ resembles the $N$-soliton formula, which also can be written in the form \eqref{eq:tau-form} with a matrix of the form $\overline{\mathbf{B}}=\mathbb{I}-\ii \mathbf{H}(t,x)/\epsilon$; here $\mathbf{H}(t,x)$ is a Hermitian matrix having simple structure (see \cite[Section 3.1]{MatsunoBook} and \cite[Eqn.\@ 1.8]{BGM23}). It is also true that the $N$-soliton formula is a rational function of $x$ when $t=0$ (see for instance~\cite[Eqn.\@ 1.19]{Sun2020}).  However, the $N$-soliton solutions remain rational in $x$ for all time $t>0$, while in general this is not true of the solution $u(t,x)$ of the Cauchy problem with rational initial data.  We wish to emphasize that the exact formul\ae\ for $u(t,x)$ given in Theorem~\ref{thm:inversion-formula} represent also solutions that in the setting of the Fokas-Ablowitz inverse-scattering transform~\cite{FokasAblowitz83,KaupMatsuno98} correspond to nonzero reflection coefficients.  See \cite{MillerW16a} where the scattering data are computed explicitly for general rational initial conditions of the form \eqref{eq:rationalIC}; generally these solutions consist of both reflectionless/solitonic and reflective/dispersive components.  This makes the explicit formula~\eqref{eq:lambda-formula} especially remarkable.
    \label{rem:remarkable}
\end{remark}

\begin{remark}[An explicit invariant space for $u(t,x)$]
As pointed out in Remark~\ref{rem:remarkable} above, the linear subspace $L^2(\mathbb{R})$ consisting of rational functions with simple poles is not invariant under the Benjamin-Ono flow.  However, one may define explicitly a \emph{nonlinear} subset of $L^2(\mathbb{R})$ forward-invariant under \eqref{eq:BO} simply as the set of functions of $x$ given by \eqref{eq:u-Piu}--\eqref{eq:lambda-formula} for any fixed $t>0$ and arbitrary $u_0$ of the form \eqref{eq:rationalIC}.  According to Remark~\ref{rem:tzero-and-negative}, this set clearly includes all rational functions of the form \eqref{eq:rationalIC} (by taking $t\downarrow 0$) and can be further made backward-invariant by replacing $t$ with $-t$ as indicated.
    \label{rem:invariant space}
\end{remark}

\begin{remark}[Implementation]
Since the elements of the matrices $\mathbf{A}(t,x)$ and $\mathbf{B}(t,x)$ involve the specific poles $\{p_n\}_{n=1}^N$ of $u_0$ in the upper half-plane, if one is given $u_0$ as a rational function in the form $u_0(x)=P(x)/Q(x)$ with $\deg(P)<\deg(Q)$ to guarantee $u_0\in L^2(\mathbb{R})$, one must first have access to the factorization of $Q(x)$ into simple factors in order to be able to use Theorem~\ref{thm:inversion-formula}.
    \label{rem:FundamentalTheoremOfAlgebra}
\end{remark}

Next we study the solution for $\epsilon=1$ generated by minus a soliton, i.e., with initial data   
\begin{equation}u_0(x)=\frac{-2}{1+x^2}\ .
\label{eq:minus-soliton}
\end{equation}
To formulate our result, we first introduce a certain spectral transform of general $u_0$ with $\langle x\rangle u_0\in L^2(\mathbb{R})$ associated with the Lax operator (here written for the case $\epsilon=1$) $L_{u_0}:=-\ii \partial_x -\Pi u_0$.  For given $\lambda>0$, one first determines the unique solution $m=m_-(x,\lambda)$ of $L_{u_0}m=\lambda m$ that is the boundary value for $x\in\mathbb{R}$ of a function analytic and bounded in the upper half $x$-plane and that has the asymptotic behavior $\ee^{-\ii\lambda x}m_-(x,\lambda)\to 1$ as $x\to -\infty$.  Then a \emph{distorted Fourier transform} may be defined on the image of $L^2(\mathbb{R})$ under $\Pi$ by 
\begin{equation}
    f\mapsto \widetilde{f}(\lambda):=\int_\mathbb{R}f(x)m_-(x,\lambda)^*\,\dd x,\quad\lambda>0.
\end{equation}
Note that in the case $u_0\equiv 0$, we have $m_-(x,\lambda)=\ee^{\ii\lambda x}$ and hence $\widetilde{f}(\lambda)=\widehat{f}(\lambda)$ reduces to the classical Fourier transform. All the most important properties of the map $f\mapsto\widetilde{f}$ needed for our analysis are described below in Section~\ref{section:spectral-Lax}.  Then a spectral function $\alpha\in L^2(0,\infty)$ is defined as follows:  
\begin{equation}
    \alpha(\lambda):=\widetilde{\Pi u_0}(\lambda),\quad\lambda>0.
\end{equation}

As an application of the exact solution formula given in Theorem~\ref{thm:inversion-formula}, we will prove the following result.  Here $\ee^{t\partial_x |D_x|}$ denotes the unitary map taking given initial data in $L^2(\mathbb{R})$ to the corresponding solution of the linearized Benjamin-Ono equation for $\epsilon=1$ given by $u_t = \partial_x |D_x|u$.
\begin{theorem}[Long-time asymptotics with initial datum equal to minus a soliton]\label{thm:minussoliton}
Let $u$ be the solution of \eqref{eq:BO} with $\epsilon =1$ and initial data \eqref{eq:minus-soliton}.
Then we have
\begin{equation}\lim_{t\to+\infty}\left\|u(t,\diamond)-{\rm e}^{t\partial_x|D_x|}u_+(\diamond)\right\|_{L^2(\mathbb{R})}= 0\ ,
\label{eq:L2-convergence}
\end{equation}
where $u_+$ is the real-valued $L^2(\mathbb{R})$ function characterized by its Fourier transform as follows:
\begin{equation}
    \forall \xi >0\ ,\ \widehat{ u_+}(\xi )=\alpha(\xi).
\end{equation}
\end{theorem}

\begin{corollary}[Scattering in $H^s(\mathbb{R})$ for any $s\ge 0$]
The convergence \eqref{eq:L2-convergence} also holds in $H^s(\mathbb{R})$ for any $s\ge 0$.
\end{corollary}
\begin{proof}
    For any $f\in H^s(\mathbb{R})$ for all $s\ge 0$, we have $\|f\|_{H^s(\mathbb{R})}\le \|f\|_{H^{2s}(\mathbb{R})}^{1/2}\|f\|_{L^2(\mathbb{R})}^{1/2}$ by interpolation.
    Note that with data \eqref{eq:minus-soliton}, the solution $u(t,\diamond)\in H^s(\mathbb{R})$ for all $s\ge 0$ uniformly in $t\ge 0$ because all $H^s$-norms are controlled by quantities conserved under \eqref{eq:BO}.  Also, $u_+\in H^s(\mathbb{R})$ for all $s$ because $\alpha(\xi)$ decays exponentially as $\xi\to+\infty$ (see \eqref{eq:minus-soliton-alpha} below), a property that is preserved uniformly in $t\ge 0$ by the Fourier multiplier $\ee^{t\partial_x|D_x|}$.  Therefore for given $s\ge 0$, $u(t,\diamond)$ and $\ee^{t\partial_x|D_x|}u_+(\diamond)$ are uniformly bounded in $H^{2s}(\mathbb{R})$, so taking $f(\diamond)=u(t,\diamond)-\ee^{t\partial_x|D_x|}u_+(\diamond)$ and using \eqref{eq:L2-convergence}, the result follows.
\end{proof}

\begin{remark}[Recent related literature]
 For some recent references where similar long-time asymptotic results have been obtained, we mention the work of Perry \cite{Perry16} and the $L^2$-refinement of  Nachman, Regev, and Tataru \cite{NachmanRT20} establishing scattering for the defocusing Davey-Stewartson type II equation; the works of Duyckaerts, Kenig, and Merle \cite{DuyckaertsKM23} and of Jendrej and Lawrie \cite{JendrejL23} establishing soliton resolution for the radial nonlinear wave equation; and the work of Kim and Kwon \cite{KimK24} establishing soliton resolution for the Calogero-Moser derivative nonlinear Schr\"odinger equation.   In particular, in \cite{Perry16,NachmanRT20} the analogue of $u_+$ for the defocusing Davey-Stewartson type II equation is identified in terms of a reflection coefficient associated to a $\overline{\partial}$-problem involving the initial data $u_0$, which like $\alpha(\lambda)$ can also be viewed as a distortion of the Fourier transform.  All of the aforementioned models are either mass-critical or energy-critical, whereas the Benjamin-Ono equation is noncritical with respect to all of its conserved quantities.
 Long-time dispersive estimates were previously shown to hold for the Benjamin-Ono equation in the case of small data by Ifrim and Tataru \cite{IfrimTataru2019}.  
\end{remark}

\begin{remark}[$\alpha(\lambda)$ versus the reflection coefficient $\beta(\lambda)$]
    In the Fokas-Ablowitz theory of the Benjamin-Ono equation \cite{FokasAblowitz83} (see also \cite{KaupMatsuno98,Wu16,Wu17}), one associates four distinct Jost solutions with the Lax operator $L_{u_0}$, denoted $M(x,\lambda)$, $N(x,\lambda)$, $\overline{M}(x,
\lambda)$, and $\overline{N}(x,\lambda)$ for $\lambda>0$.  Of these, $m=\overline{M}(x,\lambda)$ and $m=N(x,\lambda)$ are solutions of $L_{u_0}m=\lambda m$ bounded and analytic in the upper half $x$-plane with asymptotic behavior (for $\epsilon=1$) $\overline{M}(x,\lambda)\ee^{-\ii\lambda x}\to 1$ as $x\to -\infty$ and $N(x,\lambda)\ee^{-\ii\lambda x}\to 1$ as $x\to +\infty$.  Therefore $m_-(x,\lambda)$ coincides with $\overline{M}(x,\lambda)$, and $N(x,\lambda)$ is proportional to $\overline{M}(x,\lambda)$ \cite{KaupMatsuno98} (see also Proposition~\ref{geneigen} below).  On the other hand, $m=M(x,\lambda)$ and $m=\overline{N}(x,\lambda)$ are solutions of the modified equation $L_{u_0}m=\lambda(m-1)$ analytic and bounded in the upper half $x$-plane with asymptotic behavior $M(x,\lambda)\to 1$ as $x\to -\infty$ and $\overline{N}(x,\lambda)\to 1$ as $x\to +\infty$.  Unlike $\overline{M}$ and $N$, $M$ and $\overline{N}$ have analytic character in $\lambda$, with $M$ and $\overline{N}$ being the boundary values for $\lambda>0$ of a common meromorphic function $\lambda\mapsto W(x,\lambda)$ defined for $\lambda\in \mathbb{C}\setminus\mathbb{R}_+$ (and having poles only at the negative point spectrum of $L_{u_0}$) from the upper and lower half $\lambda$ planes, respectively.  The function $W$ has a jump across the positive real $\lambda$-axis given by the difference $M(x,\lambda)-\overline{N}(x,\lambda)$ which can be expressed for $\lambda>0$ in the form $\beta(\lambda)N(x,\lambda)$, ultimately leading to a nonlocal Riemann-Hilbert jump condition that is the centerpiece of the inverse-scattering theory from this point of view.  The spectral function $\beta(\lambda)$ is called the \emph{reflection coefficient} and it can be expressed in terms of the Jost solution $M(x,\lambda)$ as follows:
\begin{equation}
    \beta(\lambda)=\ii\int_\mathbb{R} u_0(x)M(x,\lambda)\ee^{-\ii\lambda x}\,\dd x,\quad\lambda>0.
\end{equation}
One can check that if $u_0\equiv 0$ is used to construct $M(x,\lambda)$, one gets $M(x,\lambda)\equiv 1$, and thus $-\ii\beta(\lambda)$ is also a kind of distorted Fourier transform of $u_0$.  We are however not aware of any general relationship between $\alpha(\lambda)$ and $\beta(\lambda)$.
\end{remark}

\begin{remark}[Exact formul\ae\ for $\alpha(\lambda)$ and $\beta(\lambda)$ when $u_0(x)$ is given by \eqref{eq:minus-soliton}]
In our proof of Theorem~\ref{thm:minussoliton} we first analyze $u(t,2ty)$ for $y<0$ and $t\to +\infty$ and show that 
\begin{equation}
    \ee^{-\ii\pi/4}\ee^{\ii ty^2}\sqrt{4\pi t}\Pi u(t,2ty)\to \frac{2\pi\ee^{-y}}{\mathrm{Ei}(-2y)+\ii\pi}
    \label{eq:intro-locally-uniform}
\end{equation}
holds locally uniformly in $y$ (see Proposition~\ref{prop:weak-asymptotics} in Section~\ref{section:weak-asymptotics}).  Here $\mathrm{Ei}(k)$ is the special function defined by 
\begin{equation}
    \mathrm{Ei}(k):=\mathrm{P.V.}\int_{-\infty}^k\frac{\ee^s}{s}\,\dd s,\quad k>0.
    \label{eq:Ei}
\end{equation}
See \cite[Eqn.\@ 6.2.5]{DLMF}. We also compute directly the spectral function $\alpha(\lambda)$ in Section~\ref{section:spectral-Lax} and observe that
\begin{equation}
    \alpha(\lambda)=\frac{2\pi\ee^\lambda}{\mathrm{Ei}(2\lambda)+\ii\pi},\quad\lambda>0,
    \label{eq:minus-soliton-alpha}
\end{equation}
so the limit on the right-hand side of \eqref{eq:intro-locally-uniform} could be written as $\alpha(-y)$.
The locally-uniform convergence described by \eqref{eq:intro-locally-uniform} is illustrated in Figure~\ref{fig: hatPsiConv}, where the left-hand side is computed from contour integration using Theorem~\ref{thm:inversion-formula} for various values of $t$ and the right-hand side is shown with a red curve.
\begin{figure}
\begin{center}
    \includegraphics[width=0.48\linewidth]{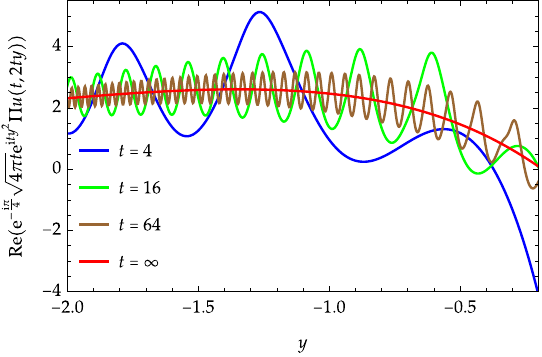}\hfill%
    \includegraphics[width=0.48\linewidth]{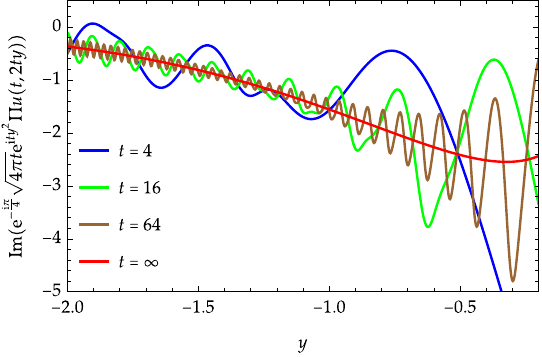}
\end{center}
\caption{A demonstration of the convergence guaranteed by Proposition~\ref{prop:weak-asymptotics}.  In the left/right panel, the red curve is the real/imaginary part of $\alpha(-y)$, respectively, see~\eqref{eq:minus-soliton-alpha}.}
\label{fig: hatPsiConv}
\end{figure}

The reflection coefficient $\beta(\lambda)$ has been computed for general initial data of the form \eqref{eq:rationalIC} in \cite{MillerW16a}.  Specializing that result for $\epsilon=1$ in the case of initial data $u_0(x)$ given by \eqref{eq:minus-soliton}, one finds that
\begin{equation}
    \beta(\lambda)=\frac{2\pi\ii\ee^\lambda}{\mathrm{Ei}(2\lambda)-\ii\pi},\quad\lambda>0.
\end{equation}
The result of Theorem~\ref{thm:minussoliton} could therefore have been expressed explicitly in terms of either $\mathrm{Ei}(2\lambda)$ or $\beta(\lambda)$.  However, for more general initial data of the form \eqref{eq:rationalIC}, $\beta(\lambda)$ and $\alpha(\lambda)$ do not seem to be easily related; indeed note that for the data \eqref{eq:minus-soliton} we have $-\ii\beta(\lambda)=\alpha(\lambda)^*$ while for small $u_0$ we expect $-\ii\beta(\lambda)\approx \alpha(\lambda)\approx \widehat{u_0}(\lambda)$ for $\lambda>0$ instead as both are distortions of the Fourier transform of $u_0$.  Furthermore, even though following \cite{MillerW16a} one can express $\alpha(\lambda)$ for general rational initial data \eqref{eq:rationalIC} explicitly in terms of determinants of contour integrals (see \cite{LongTime-BO} for details), is it not generally possible to express the contour integrals that appear in terms of known special functions.  The reason that we choose to express Theorem~\ref{thm:minussoliton} in terms of $\alpha(\lambda)$ is that our work in progress \cite{LongTime-BO} shows that the result is true in this form more generally.  Thus we regard Theorem \ref{thm:minussoliton} as the first step towards the proof of the soliton resolution conjecture for the Benjamin-Ono equation.  Also, the fact that $\alpha(\lambda)$ appears in general instead of $\beta(\lambda)$ suggests that perhaps a second Riemann-Hilbert representation of $u(t,x)$ could be found that takes $\alpha(\lambda)$ as data rather than $\beta(\lambda)$.
\end{remark}

The rest of the paper consists of the proofs of Theorem~\ref{thm:inversion-formula} and Theorem~\ref{thm:minussoliton}. The proof of Theorem~\ref{thm:inversion-formula} relies on a more generally-applicable solution formula for the Benjamin-Ono equation found by one of the authors in~\cite{Gerard22}. We recall this formula in Section~\ref{section:inversion}. Then in Section~\ref{section:inversion-rational}, we invert the operator involved in the solution formula when the initial data is rational, and we find conditions for the inverse to belong to the Hardy space in Section~\ref{section:holomorphy-conditions}. We conclude that formul\ae~\eqref{eq:lambda-formula} and~\eqref{eq:tau-form}  hold and we show the non-vanishing of $\det(\mathbf{B})$ and $\det(\overline{\mathbf{B}})$ in  Section~\ref{section:thm-proof}.  

The proof of Theorem~\ref{thm:minussoliton} relies on two steps. First, applying Theorem~\ref{thm:inversion-formula} to the initial data \eqref{eq:minus-soliton} with $\epsilon=1$, we establish long time asymptotics in the sense of weak convergence in $L^2(-\infty, 0)$ in Section~\ref{section:weak-asymptotics}. Then, in Section~\ref{section:spectral-Lax} we develop a spectral theory for the Lax operator associated to any function $u_0$ with $\langle x\rangle u_0\in L^2(\R )$ and observe that the weak limit obtained in the first step is connected to the spectral function $\alpha(\lambda)$ as the distorted Fourier transform of $\Pi u_0$, which therefore has the same $L^2$ norm as the initial datum.

In Appendix~\ref{section:soliton} we show how Theorem~\ref{thm:inversion-formula}, when applied to the initial datum $u_0(x)=2/(x^2+1)$ (essentially the general case for $N=1$ with negative imaginary $c_1$), reproduces the well-known soliton solution in the case that also $\epsilon=1$, but yields nontrivial results for other values of $\epsilon>0$.  In particular our method reveals what appears to be a new explicit formula for the $M$-soliton solution of the Benjamin-Ono equation in terms of just two contour integrals (see \eqref{eq:u-N1-generalM}). 

\subsection*{Acknowledgements}
The authors acknowledge their sponsors:
\begin{itemize}
\item E. Blackstone was partially supported by the National Science Foundation under grant DMS-1812625.
\item L. Gassot was supported by the France 2030 framework program, the Centre Henri Lebesgue ANR-11-LABX-0020-01, and the ANR project HEAD--ANR-24-CE40-3260.
\item P. G\'erard was partially supported by the French Agence Nationale de la Recherche under the ANR project ISAAC--ANR-23--CE40-0015-01.
\item
P. D. Miller was partially supported by the National Science Foundation under grant DMS-2204896.  
\end{itemize}
The authors would also like to thank the Isaac Newton Institute for Mathematical Sciences, Cambridge, for support and hospitality during the program ``Emergent phenomena in nonlinear dispersive waves'', where some work on this paper was undertaken. This program was supported by EPSRC grant EP/R014604/1.

\section{General solution formula}\label{section:inversion}

In this section, we recall the solution formula for equation~\eqref{eq:BO} derived in~\cite{Gerard22}. 

The operator $\Pi$ denotes the orthogonal projection from $L^2(\mathbb{R})$ to the Hardy space $L^2_+(\mathbb{R})$ of $L^2(\mathbb{R})$ functions with Fourier transform vanishing for frequencies $\xi<0$.  Hence in particular since $u(t,x)$ is real-valued, 
    $u(t,x)=2\mathrm{Re}(\Pi u(t,x))$, and for the rational initial condition~\eqref{eq:rationalIC},
\begin{equation}
\Pi u_0(y)=\sum_{n=1}^N\frac{c_n^*}{y-p_n^*}.
\label{eq:Piu0}
\end{equation}
The Hardy space $L^2_+(\R)$ identifies with the holomorphic  functions of the complex upper half-plane $\C_+$ satisfying
\begin{equation}
\sup_{v>0}\int_{\mathbb{R}}|f(u+\ii v)|^2\,\dd u<\infty,
\label{eq:HardyNorm}
\end{equation}
in which case the Szeg\H{o} projection $\Pi$ identifies with the orthogonal projector onto this set of holomorphic functions.

For the initial-value problem written in the form \eqref{eq:BO} with initial data $u_0\in L^2(\mathbb{R})$, the solution  formula from \cite{Gerard22} reads: for every $x\in \C_+$,
\begin{equation}
\Pi u(t,x)=\frac{1}{2\ii\pi}I_+[(X^*+2\ii t\varepsilon \partial_y + 2t T_{u_0}-x \operatorname{Id})^{-1}(\Pi u_0)].
\label{eq:Gerard}
\end{equation}
Here the ingredients are the following:
\begin{itemize}
    \item $I_+$ is the linear functional defined on the subspace of $L^2(\mathbb{R})$ consisting of functions $f$ with Fourier transforms $\widehat{f}$ right-continuous at zero frequency by $I_+[f]:=\widehat{f}(0+)$.  
    \item $x\operatorname{Id}$ is a constant multiple of the identity because on the right-hand side of this formula $x\in\mathbb{C}_+$ is simply a parameter.  (On the left-hand side however, $x$ denotes the $L^2_+(\mathbb{R})$ independent variable with respect to which we understand the projection $\Pi$.)
    \item $T_{u_0}$ denotes the Toeplitz operator on $L^2_+(\mathbb{R})$ associated with multiplication by $u_0$:  
    \begin{equation}
    T_{u_0}f := \Pi(u_0f)
    \label{eq:Toeplitz}
    \end{equation}
    for $u_0\in L^\infty(\mathbb{R})$ and $f\in L^2_+(\mathbb{R})$.
    \item $\partial_y$ denotes the differential operator with respect to the independent variable $y$ of functions in the Hardy space $L^2_+(\mathbb{R})$.
    \item $X^*$ denotes the $L^2_+(\mathbb{R})$ adjoint of the multiplication $X$ of a Hardy function $f(y)$ by $y$.  The functional $I_+$ is defined on the domain of $X^*$, and the action of $X^*$ is explicitly given by~\cite[page 7]{Gerard23}
\begin{equation}
X^* f(y):=yf(y)+\frac{1}{2\ii\pi}I_+[f].
\label{eq:G-def}
\end{equation}
The second term here is a constant function; the assertion that $f\in L^2_+(\mathbb{R})$ is in the domain of $X^*$ is equivalent to the assertion that there exists a constant $\const\in\mathbb{C}$ such that $yf(y)-\const$ is in the Hardy space $L^2_+(\mathbb{R})$.
\end{itemize}
Note that the operator $X^*+2\ii t\epsilon\partial_y + 2t T_{u_0}-x\mathrm{Id}$ can be written as $X^*-x\mathrm{Id}-2t L_{u_0}$ where $L_{u_0}=-\ii\epsilon\partial_y -T_{u_0}$ denotes the Lax operator for the Benjamin-Ono equation.

\section{The action of the resolvent \texorpdfstring{of $X^*+2\ii t\epsilon\partial_y+2tT_{u_0}$}{} for rational \texorpdfstring{$u_0$}{u0}}\label{section:inversion-rational}

The main task in implementing the formula \eqref{eq:Gerard} is to invert the operator $X^*+2\ii t\epsilon\partial_y + 2tT_{u_0}-x\operatorname{Id}$ on the Hardy function $\Pi u_0(y)$. In other words, we seek a function $y\mapsto f(y;t,x)$ in the Hardy space $L^2_+(\mathbb{R})$ and also in the domain of $\partial_y$ and $X^*$ such that
\begin{equation}
(X^*+2\ii t\varepsilon \partial_y + 2t T_{u_0}-x \operatorname{Id})f=\Pi u_0.
\label{eq:eqn-for-f}
\end{equation}
We will frequently omit the parameters $t,x$ and simply write $f=f(y)$.

In this section, we prove the following lemma.
\begin{lemma} The solution $f$ to~\eqref{eq:eqn-for-f} has the following integral representation:
\begin{equation}
f(y)
	=-\frac{\ii}{2t\varepsilon} \ee^{\ii h(y)/\varepsilon}\int_{\ee^{{3\ii\pi/4}}\infty}^y  \left( u_0(z)+\const+\sum_{n=1}^N\frac{V_n}{z-p_n}\right)\ee^{-\ii h(z)/\varepsilon}\dd z
 \label{eq:f-formula}
\end{equation}
with
\begin{equation}
    \const=-\frac{1}{2\ii\pi}I_+[f],\quad V_n=2tc_nf(p_n)-c_n,\quad n=1,\dots,N.
\label{eq:lambda-V}
\end{equation}
Here, the path of integration originates at infinity in the indicated direction to the left of the branch cuts $\Gamma_1,\dots,\Gamma_N$ of $h$, and the path lies in the domain of analyticity of $h(z)$. Then $f$ has the same domain of analyticity as $h$, namely $y\in\mathbb{C}\setminus (\Gamma_1\cup\cdots\cup \Gamma_N\cup \bar{\Gamma}_1\cup\cdots\cup\bar{\Gamma}_N)$, as shown in the left-hand panel of Figure~\ref{fig:BO-BasicContours}.
\end{lemma}

\begin{proof}
If $f(\diamond)\in L_+^2(\mathbb{R})$, one can see that for $p\in\C_+$ and $y\in\R$, we have
\begin{equation}
\frac{f(y)}{y-p}=\frac{f(y)-f(p)}{y-p}+\frac{f(p)}{y-p},
\end{equation}
where the first term on the right-hand side has a holomorphic extension to the upper half-plane  $\C_+$ given by $z\in\C_+\mapsto (f(z)-f(p))/(z-p)$, and the second term in the right-hand side has a holomorphic extension on $\mathbb{C}_-$ given by  $z\in \C_-\mapsto f(p)/(z-p)$. Hence the orthogonal projection of  $y\in\R\mapsto f(y)/(y-p)$ on the Hardy space $L_+^2(\mathbb{R})$ is simply the first term:
\begin{equation}
\Pi\left(\frac{f(\diamond)}{\diamond-p}\right)(y)
	=\frac{f(y)-f(p)}{y-p}.
\end{equation}
By linearity, we sum over the partial-fraction expansion \eqref{eq:rationalIC} of $u_0(y)$ and hence infer that
\begin{equation}
T_{u_0}f(y)=u_0(y)f(y)-\sum_{n=1}^N\frac{c_nf(p_n)}{y-p_n}.
\end{equation}
Using this along with \eqref{eq:Piu0} and \eqref{eq:G-def} in \eqref{eq:eqn-for-f}, and writing $\Pi u_0(y)=u_0(y)-(\operatorname{Id}-\Pi)u_0(y)$ on the right-hand side gives: for $y\in\C_+$,
\begin{equation}
2\ii t\epsilon f'(y) + \left[y-x+2tu_0(y)\right]f(y)=u_0(y)+\const+\sum_{n=1}^N\frac{V_n}{y-p_n},
\label{eq:f-ODE}
\end{equation}
in which the constants $\const$ and $V_1,\dots,V_N$ are given by~\eqref{eq:lambda-V}.

Now we observe that $\const$ and $V_1,\dots,V_N$ depend on the unknown function $f$, but they are constants and hence the form of the right-hand side of \eqref{eq:f-ODE} as a function of $y$ is explicit.  This allows us to view \eqref{eq:f-ODE} as a first-order linear differential equation for $f(y)$, the solution of which will depend on the constants $\const$ and $V_1,\dots,V_N$.  We will show in Section~\ref{section:holomorphy-conditions} that the constants are then uniquely determined so that when $\mathrm{Im}(x)\ge 0$, $f(y)$ is a function in $L_+^2(\mathbb{R})$ that is also in the domain of $X^*$ and of $\partial_y$.  In particular, since \eqref{eq:Gerard} can be written as 
\begin{equation}\label{eq:Pi-lambda}
\Pi u(t,x)=\frac{I_+[f]}{2\ii\pi}=-\const,
\end{equation}
the value of the constant $\const$ is exactly what is needed to solve the Cauchy problem \eqref{eq:BO}--\eqref{eq:rationalIC}:
\begin{equation}
    u(t,x)=2\mathrm{Re}(\Pi u(t,x)) = -2\mathrm{Re}(\const).
    \label{eq:u-from-lambda}
\end{equation}

An integrating factor for \eqref{eq:f-ODE} is $\ee^{-\ii h(z)/\epsilon}$,  where $h$ is defined for $z\in\mathbb{C}\setminus (\Gamma_1\cup\cdots\cup \Gamma_N\cup\bar{\Gamma}_1\cup\cdots\cup\bar{\Gamma}_N)$ as explained in the introduction. We assume that $\epsilon>0$ and also that $t>0$.  Then, the dominant quadratic term in $h(z)$ behaves such that the function $\ee^{-\ii h(z)/\epsilon}$ decays exponentially to zero in the directions $\arg(z)=3\pi/4$ and $\arg(z)=-\pi/4$, so the use of the integrating factor
leads to the integral representation of a particular solution in the form~\eqref{eq:f-formula}.
\end{proof}

\section{Conditions for \texorpdfstring{$f$}{f} to lie in the Hardy space}\label{section:holomorphy-conditions}

Now, assuming that $\mathrm{Im}(x)\ge 0$, we seek conditions on $\const$ and $V_1,\dots,V_N$ such that $f$ is a Hardy function of $y$.  This requires two things:  \begin{enumerate}[(i)]
\item\label{i} that $f(y)$ is analytic at each of the points $y=p_1,\dots,p_N$ (guaranteeing that $f(y)$ is analytic except on the lower half-plane branch cuts $\bar{\Gamma}_1,\dots,\bar{\Gamma}_N$, and in particular for 
$\mathrm{Im}(y)>-\delta$ for some $\delta>0$);
\item\label{ii} that the $L_+^2(\mathbb{R})$-norm~\eqref{eq:HardyNorm} of $f(y)$ is finite.
\end{enumerate}

For
condition~\ref{ii} to hold given~\ref{i}, it will be enough to prove that $yf(y)$ has a limit as $y\to\infty$ in a horizontal strip $|\mathrm{Im}(y)|<\delta$ centered on the real line. By confirming that this limit is the same in the left and right directions, we will further conclude that $f(y)$ is in the domain of $X^*$.  The  equation \eqref{eq:eqn-for-f} then implies that $f(y)$ is also in the domain of $\partial_y$.

For condition~\ref{i}, we show the following lemma. Recall the contours $C_m$ for $m=1,\dots,N$ defined in Definition~\ref{def:C}.

\begin{lemma}
    Assume the conditions 
\begin{equation}
    I_m:=\int_{C_m}\left(u_0(z)+\const+\sum_{n=1}^N\frac{V_n}{z-p_n}\right)\ee^{-\ii h(z)/\epsilon}\,\dd z = 0,\quad m=1,\dots,N.
\label{eq:analyticity}
\end{equation}
Then the formula \eqref{eq:f-formula} defines $f(y)$ as a function analytic for $\mathrm{Im}(y)>-\delta$.  
\label{lem:analyticity}
\end{lemma}

For the proof, we apply the reasoning behind \cite[Corollary 1]{MillerW16a}.
Before getting into the details, let us first give an idea of the proof in the simplest situation: $N=1$.  There are four cases to consider.
\paragraph{The exceptional case: $\ii c_1/\epsilon\in \mathbb{N}_{<0}$} In this case, the contour $C_1$ goes from $\ee^{3\ii\pi/4}\infty$ to $p_1$ and the integrand is integrable at $z=p_1$, so that by the condition $I_1=0$ one can change the starting point of the integration path in~\eqref{eq:f-formula}:
\begin{equation}
f(y)=-\frac{\ii}{2t\epsilon}\ee^{\ii h(y)/\epsilon}\int_{p_1}^y\left(u_0(z)+\const+\frac{V_1}{z-p_1}\right)\ee^{-\ii h(z)/\epsilon}\,\dd z.
\label{eq:f-formula-rewrite-N1}
\end{equation}
This defines a holomorphic function at $z=p_1$. Indeed, since the phase has the expression $\ee^{-\ii h(z)/\epsilon}=(z-p_1)^{-\ii c_m/\epsilon}(z-p_1^*)^{-\ii c_m^*/\epsilon}\ee^{-\ii (z-x)^2/(4t\epsilon)}$, we can write the integrand in the form $H_1(z) (z-p_1)^{-\ii c_m/\epsilon-1}$ with $H_1$ being holomorphic at $p_1$. Then by Taylor expansion of $H_1$ around $z=p_1$ and term-by-term integration, we see that $f(y)$ also has a convergent Taylor series around $y=p_1$.
\paragraph{The non-exceptional case with $\operatorname{Re}(\ii c_1/\epsilon)<0$}  In this case, the integrand is also integrable at $z=p_1$. Hence we can choose the contour $C_1$ to be the concatenation  of $C_1^\searrow$ from $\ee^{3\ii\pi/4}\infty$ to $p_1$ and $C_1^\nwarrow$ from $p_1$ to $\ee^{3\ii\pi/4}\infty$ along opposite sides of the branch cut $\Gamma_1$.
Correspondingly,
the integral $I_1$ decomposes into $I_1=I_1^{\searrow}+I_1^{\nwarrow}$. The logarithmic jump being constant along $\Gamma_1$, there holds $I_1=(1-\ee^{2\pi\ii(-\ii c_m/\epsilon-1)})I_1^\searrow$ so that $I_1^\searrow$ vanishes. The proof of the exceptional case then applies and we conclude that $f$ is holomorphic at $z=p_1$.
\paragraph{The non-exceptional case with $\ii c_1/\epsilon=M\in\mathbb{N}$} In this case, the integrand in~\eqref{eq:f-formula} has a pole of order at most $M+1$ at $z=p_1$. The condition $I_1=0$ then implies the vanishing of the residue of this integrand at $z=p_1$ so that the integral in the definition \eqref{eq:f-formula} of $f(y)$ defines a function analytic in the upper half-plane except for a pole of order at most $M$  at $z=p_1$, which is compensated  by zero of order $M$ of $\ee^{\ii h(y)/\epsilon}$.  Hence again $f$ is holomorphic at $z=p_1$.
\paragraph{The non-exceptional case with $\operatorname{Re}(\ii c_1/\epsilon)\geq 0$ and $\ii c_1/\epsilon\not\in\mathbb{N}$}  In this case, we can integrate by parts a sufficient number of times to replace the factor $(z-p_1)^{-\ii c_m/\epsilon-1}$ coming from $\ee^{-\ii h(z)/\epsilon}/(z-p_1)$ in the integrand by $(z-p_1)^{-\ii c_m/\epsilon+M-1}$ for $M$ large enough, after which the proof is the same as in the non-exceptional case with $\operatorname{Re}(\ii c_1/\epsilon)< 0$ with the same conclusion that $I_1=0$ implies that $f$ is analytic at $y=p_1$.

Analyticity of $f$ at just the one point $y=p_1$ is enough to guarantee analyticity for all $y\in\mathbb{C}_+$ in the $N=1$ case.  With the idea of the proof clear for $N=1$, we now proceed to the general case.
\begin{proof}[Proof of Lemma~\ref{lem:analyticity}]
It is sufficient to show that the formula \eqref{eq:f-formula} defines an analytic function in the neighborhood of each of the points $p_1,\dots,p_N$.

Selecting an index $m$, first suppose that $\mathrm{Re}(\ii c_m/\epsilon)<0$.  Then, the integrand in \eqref{eq:analyticity} is integrable at $z=p_m$, and we claim that for $y$ near $p_m$, $f(y)$ can be written in the equivalent form
\begin{equation}
f(y)=-\frac{\ii}{2t\epsilon}\ee^{\ii h(y)/\epsilon}\int_{p_m}^y\left(u_0(z)+\const+\sum_{n=1}^N\frac{V_n}{z-p_n}\right)\ee^{-\ii h(z)/\epsilon}\,\dd z,
\label{eq:f-formula-rewrite}
\end{equation}
where the path of integration avoids the branch cut $\Gamma_m$ of $h(z)$ emanating from $p_m$.

Indeed, in the special case that $m$ is an exceptional index, deforming the path in \eqref{eq:f-formula} to pass through $z=p_m$ and using the condition $I_m=0$ in \eqref{eq:analyticity} leads directly to \eqref{eq:f-formula-rewrite}.

 If $\mathrm{Re}(\ii c_m/\epsilon)<0$ but $m$ is a non-exceptional index, we reason as follows.  If there exists a non-exceptional index $j<m$, let $j$ be the largest such index and observe that the difference $I_m-I_j$ from \eqref{eq:analyticity} is an integral $\overline{I}_m$ of the same integrand about a loop encircling in the counterclockwise sense only the branch cut $\Gamma_m$ of $h(z)$ and no other singularities of the integrand (the integrand is analytic at $z=p_k$ for $j<k<m$ even though $h(z)$ has singularities because these are all exceptional indices).  Otherwise, $\overline{I}_m=I_m$ is itself such an integral.  Since on opposite sides of the branch cut $\Gamma_m$ the integrand differs by the non-unit factor $\ee^{2\pi\ii(-\ii c_m/\epsilon-1)}$ and since the integrand is integrable at $z=p_m$, the resulting integral condition $\overline{I}_m=0$ implies that the integral $\overline{I}_m^\searrow$ inwards from $z=\infty$ along just the one side of $\Gamma_m$ vanishes (writing $\overline{I}_m=\overline{I}_m^\searrow+\overline{I}_m^\nwarrow$ we have $\overline{I}_m^\nwarrow=-\ee^{2\pi\ii (-\ii c_m/\epsilon-1)}\overline{I}_m^\searrow$, and $1-\ee^{2\pi\ii (-\ii c_m/\epsilon-1)}\neq 0$).  Taking the sum of this ``one-sided'' integral $\overline{I}_m^\searrow$ and $I_j$ from \eqref{eq:analyticity}, or zero if there is no exceptional index $j<m$, one obtains an integral condition on a contour terminating at $z=p_m$, just as in the case that $m$ is exceptional, that combines with the definition \eqref{eq:f-formula} to yield \eqref{eq:f-formula-rewrite}.

To make use of the formula for $f(y)$ rewritten in the form \eqref{eq:f-formula-rewrite}, we note that the integrand can be written in the form $H_m(z)(z-p_m)^{-\ii c_m/\epsilon-1}$ where $H_m(z)$ is holomorphic in a neighborhood of 
$z=p_m$ and the power function is cut along $\Gamma_m$, and that $\ee^{\ii h(y)/\epsilon}$ can be written in a similar form $G_m(y)(y-p_m)^{\ii c_m/\epsilon}$ with $G_m(y)$ analytic at $p_m$. Using these results in \eqref{eq:f-formula-rewrite}, expanding $H_m(z)$ in a Taylor series about $z=p_m$ and integrating term-by-term, the desired analyticity of $f(y)$ at $y=p_m$ is easily established under the operating assumption that $\mathrm{Re}(\ii c_m/\epsilon)<0$.  

Next, suppose that $\ii c_m/\epsilon$ is a nonnegative integer $M$.  Then the integrand in the conditions \eqref{eq:analyticity} has a pole of order at most $M+1$ at $z=p_m$, while the factor $\ee^{\ii h(y)/\epsilon}$ is analytic with a zero of order $M$ at $z=p_m$.  Again letting $j$ denote the largest non-exceptional index with $j<m$ as before, the difference $\overline{I}_m=I_m-I_j$ (or just $\overline{I}_m=I_m$ if there is no such index $j$) then yields an integral condition $\overline{I}_m=0$ on a contour that can be deformed to a loop surrounding the point $z=p_m$, which is the only singularity of the integrand therein.  Therefore, the residue of the meromorphic integrand vanishes at $z=p_m$ and hence the integral factor in the definition \eqref{eq:f-formula} is a single-valued meromorphic function of $y$ near $p_m$ with a pole of order $M$, and again analyticity of $f(y)$ at $y=p_m$ follows.  

Finally, we suppose that $\mathrm{Re}(\ii c_m/\epsilon)\ge 0$ but that $\ii c_m/\epsilon$ is not an integer. Then, both in the definition \eqref{eq:f-formula} and the integrals $I_k$ from \eqref{eq:analyticity} with non-exceptional indices $k$,
we integrate by parts $M$ times using the representation of the integrand
\begin{equation}
H_m(z)(z-p_m)^{-\ii c_m/\epsilon-1} = \frac{H_m(z)}
{(-\ii c_m/\epsilon)_M}\frac{\dd^M}{\dd z^M}(z-p_m)^{-\ii c_m/\epsilon+M-1}
\end{equation}
where $H_m(z)$ is analytic at $z=p_m$, $(a)_M:=\Gamma(a+M)/\Gamma(a)$ denotes Pochhammer's symbol, and $M$ is sufficiently large that $\mathrm{Re}(-\ii c_m/\epsilon)+M>0$.  No boundary terms are produced for the integrals $I_k$ because under the indicated condition on $k$, the contour $C_k$ has no finite endpoints.  In the formula \eqref{eq:f-formula} however, each integration by parts will produce a contribution from the finite endpoint $z=y$.  In this case, one uses the representation $\ee^{\ii h(y)/\epsilon}=G_m(y)(y-p_m)^{\ii c_m/\epsilon}$ with $G_m(y)$ analytic at $y=p_m$ to conclude that the contribution of each of these terms to $f(y)$ is holomorphic at $y=p_m$.  For the remaining integral contribution to $f(y)$, the integrand is integrable at $z=p_m$ and one uses exactly the same argument as in the case $\mathrm{Re}(\ii c_m/\epsilon)<0$ with the rewritten integrals $I_k$ to conclude finally that $f(y)$ is analytic at $y=p_m$.
\end{proof}

\begin{remark}
    Since Lemma~\ref{lem:analyticity} already guarantees that $f$ is analytic on the potentially-complicated branch cuts $\Gamma_1,\dots,\Gamma_N$, for the rest of the argument it is sufficient to take the branch cuts of $h(z)$ to be straight rays from each of the points $p_j$ to $\infty$ with $\arg(z)=3\pi/4$ and to define $h$ in the lower half-plane by Schwarz symmetry $h(z^*)=h(z)^*$, and we shall do so.
    \label{rem:change-cuts}
\end{remark}

To deal with condition (ii), we let $C_0$ be as in Definition~\ref{def:C}.  See Figure~\ref{fig:BO-BasicContours}, right-hand panel.  Then we have the following.  
\begin{lemma}
    Let $\mathrm{Im}(x)\ge 0$.  Suppose that the conditions \eqref{eq:analyticity} hold, and also that
\begin{equation}
I_0:=\int_{C_0}\left(u_0(z)+\const +\sum_{n=1}^N\frac{V_n}{z-p_n}\right)\ee^{-\ii h(z)/\epsilon}\,\dd z=0.
\label{eq:decay}
\end{equation}
Then $f(y)$ defined by \eqref{eq:f-formula} lies in the Hardy space $L_+^2(\mathbb{R})$, and also in the domain of $X^*$.
\label{lem:decay}
\end{lemma}

Before getting into the details, we first describe the main ideas of the proof.  We first show that $f(y)\sim \const/y$ uniformly as $y\to \infty$ in a strip of the form $0\leq\operatorname{Im}(y)<\delta$ for some $\delta>0$. We prove that this asymptotic expansion always holds as $\operatorname{Re}(y)\to-\infty$, but we will need the condition~\eqref{eq:decay} for the expansion to hold as $\operatorname{Re}(y)\to+\infty$.

We start from expression~\eqref{eq:f-formula} written in the form
\begin{equation}
f(y)
	=-\frac{\ii}{2t\varepsilon} \int_{\ee^{{3\ii\pi/4}}\infty}^y  \left( u_0(z)+\const+\sum_{n=1}^N\frac{V_n}{z-p_n}\right)\ee^{-\ii (h(z)-h(y))/\varepsilon}\dd z.
\end{equation}
Then, we show that the dominant contribution in the exponent $-\ii (h(z)-h(y))/\epsilon$ comes from the term $-\ii ((z-x)^2-(y-x)^2)/(4t\epsilon)$, whose real part is made maximal at $z=y$ on a well-chosen path from $\ee^{3\ii\pi/4}\infty$ to $y$. More precisely, to establish an asymptotic expansion, we  use a real parametrization of this path (see~\eqref{eq:schwarz} below), and write $f$ in the form
\begin{equation}
f(y)
	=\const (y-x)\int_0^{+\infty}\ee^{-|y-x|^2 s}\phi(s)\dd s,
\end{equation} with $\phi(0)=1$.
Applying Watson's lemma (Laplace's method for exponential integrals on $\mathbb{R}_+$ with exponential factor $\ee^{-ms}$ in the integrand for large $m>0$), we conclude that $f(y)\sim\const/y$ as $\operatorname{Re}(y)\to-\infty$ and $0\leq\operatorname{Im}(y)<\delta$.

The condition~\eqref{eq:decay} implies that there also holds
\begin{equation}
f(y)
	=-\frac{\ii}{2t\varepsilon} \int_{\ee^{{-\ii\pi/4}}\infty}^y  \left( u_0(z)+\const+\sum_{n=1}^N\frac{V_n}{z-p_n}\right)\ee^{-\ii (h(z)-h(y))/\varepsilon}\dd z.
\end{equation}
Consequently, a similar study implies that $f$ has the same asymptotic expansion  $f(y)\sim\const/y$ as $\operatorname{Re}(y)\to+\infty$ and $0\leq\operatorname{Im}(y)<\delta$.

This decay property of $f$ at infinity implies that $f\in L_+^2(\mathbb{R})$, and the common limit of $yf(y)$ at $\pm\infty$ obtained from~\eqref{eq:decay} also implies that $f$ belongs to the domain of $X^*$. With this summary of the argument in hand, we proceed to the details.

\begin{proof}[Proof of Lemma~\ref{lem:decay}]
Applying Lemma~\ref{lem:analyticity}, to prove that $f\in L^2_+(\mathbb{R})$, it is enough to show that $f(y)=O(y^{-1})$ uniformly as $y\to\infty$ in a horizontal strip of the form $0\le \mathrm{Im}(y)<\delta$ for some $\delta>0$.  Indeed, in taking the supremum over $v>0$ in \eqref{eq:HardyNorm}, one can assume also that $0<v<\delta$, and then by analyticity, and hence also continuity, of $f$ on the closed strip $0\le\mathrm{Im}(y)\le\delta$, for every $L>0$ there is a constant $K^{(1)}_L>0$ such that $|f(u+\ii v)|^2\le K^{(1)}_L$ holds for $|u|\le L$ and $0\le v\le\delta$. Therefore, if $f(u+\ii v)=O(u^{-1})$ holds uniformly for $0\le v\le\delta$ as $u\to\infty$, then for every $L>0$ sufficiently large, there is a constant $K^{(2)}>0$ such that $|f(u+\ii v)|^2\le K^{(2)}u^{-2}$ holds for $|u|\ge L$ and $0\le v\le\delta$.  Consequently,
\begin{equation}
    \int_\mathbb{R}|f(u+\ii v)|^2\,\dd u\le 2LK^{(1)}_L + 2K^{(2)}\int_L^{\infty}\frac{\dd u}{u^2} = 2LK^{(1)}_L +2K^{(2)}L^{-1}<\infty
\end{equation}
holds for $0\le v<\delta$ and hence $f\in L_+^2(\mathbb{R})$. 

To show that $f(u+\ii v)=O(u^{-1})$ as $u\to\infty$ holds uniformly for $v$ bounded, we first introduce a large positive parameter by $r:=|y-x|>0$ and make the substitutions $y=x+rY$ and, in the integrand in \eqref{eq:f-formula}, $z=x+rZ$. Here, $Y$ and $Z$ are variables corresponding to $y$ and $z$ respectively, and $|Y|=1$.  
    
Then, in light of Remark~\ref{rem:change-cuts},
\begin{equation}
\begin{split}
    \ee^{\ii h(y)/\epsilon}&=\ee^{\ii (y-x)^2/(4t\epsilon)}\prod_{j=1}^N(\ee^{\ii\pi/4}(y-p_j))^{\ii c_j/\epsilon}(\ee^{-\ii\pi/4}(y-p_j^*))^{\ii c_j^*/\epsilon} \\
    &=r^{\ii P}\ee^{\ii r^2Y^2/(4t\epsilon)}
    \prod_{j=1}^N\left(\ee^{\ii\pi/4}\left(Y +\frac{x-p_j}{r}\right)\right)^{\ii c_j/\epsilon}\left(\ee^{-\ii\pi/4}\left(Y+\frac{x-p_j^*}{r}\right)\right)^{\ii c_j^*/\epsilon},
\end{split}
\end{equation}
where
\begin{equation}
    P:=\frac{1}{\epsilon}\sum_{j=1}^N (c_j+c_j^*).
    \label{eq:P-def}
\end{equation}
Likewise, the factors in the integrand in \eqref{eq:f-formula} are
\begin{equation}
    \ee^{-\ii h(z)/\epsilon}=r^{-\ii P}\ee^{-\ii r^2Z^2/(4t\epsilon)}\prod_{j=1}^N\left(\ee^{\ii\pi/4}\left(Z+\frac{x-p_j}{r}\right)\right)^{-\ii c_j/\epsilon}\left(\ee^{-\ii\pi/4}\left(Z+\frac{x-p_j^*}{r}\right)\right)^{-\ii c_j^*/\epsilon}
\end{equation}
and
\begin{equation}
\begin{split}
    u_0(z)+\const+\sum_{n=1}^N\frac{V_n}{z-p_n}&=\const + \sum_{n=1}^N \frac{c_n+V_n}{z-p_n} +\sum_{n=1}^N\frac{c_n^*}{z-p_n^*}\\
    &=\const+\frac{1}{r}\sum_{n=1}^N\left[\frac{c_n+V_n}{Z+(x-p_n)/r}+\frac{c_n^*}{Z+(x-p_n^*)/r}\right].
\end{split}
\end{equation}
Then, from \eqref{eq:f-formula}, we have
\begin{multline}
    f(y)=-\frac{\ii r}{2t\epsilon}\int_{\ee^{3\ii\pi/4}\infty}^Y
    \ee^{\ii r^2 (Y^2-Z^2)/(4t\epsilon)}\\
    {}\cdot\prod_{j=1}^N\frac{(\ee^{\ii\pi/4}(Y+(x-p_j)/r))^{\ii c_j/\epsilon}(\ee^{-\ii\pi/4}(Y+(x-p_j^*)/r))^{\ii c_j^*/\epsilon}}
    {(\ee^{\ii\pi/4}(Z+(x-p_j)/r))^{\ii c_j/\epsilon}(\ee^{-\ii\pi/4}(Z+(x-p_j^*)/r))^{\ii c_j^*/\epsilon}} \\
{}    \cdot\left(\const + \frac{1}{r}\sum_{n=1}^N\left[\frac{c_n+V_n}{Z+(x-p_n)/r} +\frac{c_n^*}{Z+(x-p_n^*)/r}\right]\right)\,\dd Z.
\label{eq:f-rescale}
\end{multline}
Here the path of integration is that obtained from the original path in the $z$-plane under the affine mapping $Z=(z-x)/r=(z-x)/|y-x|$.  The branch cuts of the factors in the denominator on the second line of \eqref{eq:f-rescale} are also mapped into the $Z$-plane from those of $h(z)$ by this scaling (a similar remark applies to the factors in the numerator as functions of $Y=(y-x)/r=(y-x)/|y-x|$).  

Now, we justify a change of variable of the form $\ii(Y^2-Z^2)/(4t\epsilon)=-s$ for $s\in\mathbb{R}_+$ by deforming the path of integration. Taking $y$ with a large negative real part and bounded imaginary part, we find that $r\gg 1$ while $Y\approx -1$, and we wish to deform the path of integration to lie along the branch of the hyperbola $\mathrm{Re}(Z^2)=\mathrm{Re}(Y^2)\approx 1$ that asymptotes to the ray $\arg(Z)=3\pi/4$.  Depending on the value of $x$ with $\mathrm{Im}(x)\ge 0$ relative to the points $p_1,\dots,p_N$, the target hyperbola may cross one or more of the images in the $Z$-plane of the straight-line branch cuts of $h(z)$ emanating from the points $p_1,\dots,p_N$, in which case we can slightly shift them to the right far from the origin in the $z$-plane to ensure analyticity of the integrand as $Z\to\infty$ along the hyperbola in the indicated direction.  Having selected the steepest descent path for which the exponential factor $\ee^{\ii r^2(Y^2-Z^2)/(4t\epsilon)}$ is positive for $Z$ along the path and maximized at the finite endpoint $Z=Y$, Laplace's method applies to give a complete asymptotic expansion of $f(y)$ in descending integer powers of $r$.  Indeed, the substitution $Z=Z(s):=Y\sqrt{1-4\ii t\epsilon s/Y^2}$, where the square root is positive when $s=0$ and analytic for $s\in\mathbb{R}$ when $Y^2\approx 1$, parametrizes the desired branch of the hyperbola by $s>0$.  Using this explicit parametrization in \eqref{eq:f-rescale} gives
\begin{multline}
    f(y)=\frac{\ii r}{2t\epsilon}\int_0^{+\infty}\ee^{-r^2s}
        {}\cdot\prod_{j=1}^N\frac{(\ee^{\ii\pi/4}(Y+(x-p_j)/r))^{\ii c_j/\epsilon}(\ee^{-\ii\pi/4}(Y+(x-p_j^*)/r))^{\ii c_j^*/\epsilon}}{(\ee^{\ii\pi/4}(Z(s)+(x-p_j)/r))^{\ii c_j/\epsilon}(\ee^{-\ii\pi/4}(Z(s)+(x-p_j^*)/r))^{\ii c_j^*/\epsilon}}\\
        {}\cdot \left(\const+\frac{1}{r}\sum_{n=1}^N\left[\frac{c_n+V_n}{Z(s)+(x-p_n)/r}+\frac{c_n^*}{Z(s)+(x-p_n^*)/r}\right]\right)Z'(s)\,\dd s.
    \label{eq:schwarz}
\end{multline}
Applying Watson's Lemma and using $Z'(0)=-2\ii t\epsilon/Y$ yields
\begin{equation}
    f(y)=\frac{\const}{rY}+O(r^{-2}) = \frac{\const}{y-x}+O(|y-x|^{-2}) = \frac{\const}{y}+O(y^{-2}),\quad y\to\infty
    \label{eq:f-large-y}
\end{equation}
with the expansion being valid for large negative $\mathrm{Re}(y)$ and bounded $\mathrm{Im}(y)$.  This result does not require the condition \eqref{eq:decay}.

To approximate $f(y)$ as $\mathrm{Re}(y)\to+\infty$ for $\mathrm{Im}(y)$ bounded, we first use \eqref{eq:decay} to write \eqref{eq:f-formula} in the form
\begin{equation}
    f(y)=-\frac{\ii}{2t\epsilon}\ee^{\ii h(y)/\epsilon}\int_{\ee^{-\ii\pi/4}\infty}^y\left(u_0(z)+\const+\sum_{n=1}^N\frac{V_n}{z-p_n}\right)\ee^{-\ii h(z)/\epsilon}\,\dd z,
    \label{eq:f-formula-southeast}
\end{equation}
where the path of integration originates from $z=\infty$ in the indicated direction to the right of all branch cuts of $h(z)$ in the lower half-plane and lies in the domain of analyticity of the integrand.  Again setting $r:=|y-x|$, rescaling by $y=x+rY$, $z=x+rZ$, taking a hyperbolic path of integration parametrized by the same function $Z=Z(s):=Y\sqrt{1-4\ii t\epsilon s/Y^2}$ for $s>0$ (although in the present situation that $\mathrm{Re}(y)$ is large and positive, $Y\approx +1$ instead of $-1$) and possibly shifting some branch cuts of $h(z)$ in the lower half-plane to the left to maintain analyticity of the integrand  along the integration path, we arrive at the same result:  the expansion \eqref{eq:f-large-y} is also valid for large positive $\mathrm{Re}(y)$ and bounded $\mathrm{Im}(y)$.  

The validity of \eqref{eq:f-large-y} for large $|\mathrm{Re}(y)|$ and bounded $\mathrm{Im}(y)$ implies that $f(u+\ii v)=O(u^{-1})$ as $u\to\pm\infty$ uniformly for bounded $v$, and hence that $f\in L_+^2(\mathbb{R})$.  However from \eqref{eq:f-large-y} we learn more, namely that $yf(y)-\const\in L_+^2(\mathbb{R})$, which shows that $f$ lies in the domain of $X^*$ as well.
\end{proof}

\section{Derivation of the solution formul\ae~\texorpdfstring{\eqref{eq:lambda-formula} and~\eqref{eq:tau-form}}{}}\label{section:thm-proof}

\subsection{Derivation of formula~\texorpdfstring{\eqref{eq:lambda-formula}}{}}
Note that equations~\eqref{eq:analyticity} and~\eqref{eq:decay} amount to a linear system of $N+1$ equations on the $N+1$ unknowns $\const, V_1,\dots,V_N$:
\begin{equation}
\mathbf{Bx}=-\mathbf{b}, 
\label{eq:linear-system}
\end{equation}
in which 
\begin{equation}
\begin{split}
\mathbf{x}&:=(\const,V_1,\dots,V_N)^\top,\\
\mathbf{b}&:=(A_{1,1},\dots, A_{N+1,1})^\top
\end{split}
\end{equation}
and $\mathbf{A}$ and $\mathbf{B}$ are defined in~\eqref{eq:Amatrix} and~\eqref{eq:Bmatrix}.
In particular, by Cramer's rule, the quantity $\const=\const(t,x)$ is given by
\begin{equation}
\const(t,x)=-\frac{\det(\mathbf{A}(t,x))}{\det(\mathbf{B}(t,x))}
\label{eq:lambda-formula-2}
\end{equation}
provided the denominator is nonzero, where $\mathbf{A}(t,x)$ is the same as $\mathbf{B}(t,x)$ but with the first column replaced by $\mathbf{b}(t,x)$. 
Note that formula~\eqref{eq:lambda-formula-2} extends by continuity to the case $\operatorname{Im}(x)=0$.
That the formula \eqref{eq:lambda-formula} solves the Cauchy problem for the Benjamin-Ono equation with rational initial data $u_0$ provided that $\det(\mathbf{B}(t,x))\neq 0$ then follows from the fact that $\Pi u=-\const$ according to~\eqref{eq:Pi-lambda}.

\subsection{Derivation of formula~\texorpdfstring{\eqref{eq:tau-form}}{}}
Next we prove the alternate formula \eqref{eq:tau-form}.
Let $j=1,\dots,N+1$.  Note that using \eqref{eq:h-def}--\eqref{eq:hprime} and \eqref{eq:Amatrix}--\eqref{eq:Bmatrix} we can obtain
\begin{equation}
\begin{split}
    \epsilon\frac{\partial B_{j1}}{\partial x}&=\int_{C_{j-1}}\left[-\ii\frac{x-z}{2t}\right]\ee^{-\ii h(z)/\epsilon}\,\dd z \\ &= \int_{C_{j-1}}\left[-\ii(u_0(z)-h'(z))\right]\ee^{-\ii h(z)/\epsilon}\,\dd z \\ &=-\ii\int_{C_{j-1}}u_0(z)\ee^{-\ii h(z)/\epsilon}\,\dd z\\ &=-\ii A_{j1},
\end{split}
\end{equation}
because $\ee^{-\ii h(z)/\epsilon}$ vanishes at each endpoint (finite or infinite) of each contour $C_{j-1}$.  An even simpler calculation shows that for $k=2,\dots,N+1$, 
\begin{equation}
    \begin{split}
        \epsilon\frac{\partial B_{jk}}{\partial x}&=\frac{\ii}{2t}\int_{C_{j-1}}\frac{z-x}{z-p_{k-1}}\ee^{-\ii h(z)/\epsilon}\,\dd z\\
        &=\frac{\ii}{2t}B_{j1}+\frac{\ii (p_{k-1}-x)}{2t}B_{jk}.
    \end{split}
\end{equation}
Therefore, by writing the derivative of the determinant as a sum over $k=1,\dots,N+1$ of the determinant obtained from $\mathbf{B}$ by replacing column $k$ by its derivative, and using the fact that all columns of $\mathbf{A}$ and $\mathbf{B}$ are the same except the first, we get
\begin{equation}
    \epsilon\frac{\partial}{\partial x}\det(\mathbf{B})=-\ii\det(\mathbf{A})+\frac{\ii}{2t}\left(\sum_{k=2}^{N+1}(p_{k-1}-x)\right)\det(\mathbf{B}),
\end{equation}
from which it follows that for any choice of the complex logarithm, 
\begin{equation}
\begin{split}
    \frac{\det(\mathbf{A})}{\det(\mathbf{B})}
    &=\ii\epsilon\frac{\partial}{\partial x}\log(\det(\mathbf{B}))+\frac{1}{2t}\sum_{k=2}^{N+1}(p_{k-1}-x)\\
    &=\ii\epsilon\frac{\partial}{\partial x}\log(\det(\mathbf{B}))+\ii\epsilon\frac{\partial}{\partial x}\log\left(\prod_{k=2}^{N+1}\ee^{\ii(x-p_{k-1})^2/(4t\epsilon)}\right).
\end{split}
\end{equation}
Hence, we obtain
\begin{equation}\label{eq:B-to-Btilde}
 \frac{\det(\mathbf{A})}{\det(\mathbf{B})}      
    =\ii\epsilon\frac{\partial}{\partial x}\log(\det(\overline{\mathbf{B}})),
\end{equation}
where the first columns of $\mathbf{B}$ and $\overline{\mathbf{B}}$ agree while the remaining columns of $\overline{\mathbf{B}}(t,x)$ are given by \eqref{eq:Btilde}. Therefore, from \eqref{eq:lambda-formula} we arrive at the representation \eqref{eq:tau-form}.

\subsection{Nonvanishing of the denominators}\label{section:nonvanishing}

The last step is to prove that $\det(\mathbf{B}(t,x))$ and $\det(\overline{\mathbf{B}}(t,x))$ are nonvanishing.  Since $\det(\mathbf{B}(t,x))$ and $\det(\overline{\mathbf{B}}(t,x))$ are explicitly related, it is enough to consider $\det(\mathbf{B}(t,x))$.  We first observe the following.
\begin{lemma}
Let $t,\epsilon>0$ be fixed.  Then the function $x\mapsto\det(\mathbf{B}(t,x))$ is entire and not identically vanishing.
\label{lem:isolated-zeros}
\end{lemma}
\begin{proof}
The contour integrals in the matrix elements of $\mathbf{B}(t,x)$ are all absolutely convergent for $t,\epsilon>0$ and arbitrary $x\in\mathbb{C}$, and the integrands are all entire functions of $x$, so it follows that $x\mapsto\det(\mathbf{B}(t,x))$ is entire.  To show that this function does not vanish identically in $x$, we write $x=R\ee^{-\ii\pi/4}$ and consider the asymptotic behavior of the matrix elements of $\mathbf{B}(t,x)$ as $R\to+\infty$.  For the first row of $\mathbf{B}(t,x)$, we observe that the factor $\ee^{-\ii (z-x)^2/(4t\epsilon)}$ in $\ee^{-\ii h(z)/\epsilon}$ has a simple critical point at $z=x=\ee^{-\ii\pi/4}R$ that may be taken to lie on the contour $C_0$ in the distant fourth quadrant of the complex $z$-plane.  After suitable rescaling of $z-x$ and noting that the logarithms in $h(z)$ agree with the principal branch when $z\approx x$, the method of steepest descent applies to show that the first row of $\mathbf{B}(t,x)$ has the form 
\begin{equation}
    \mathbf{e}_1^\top\mathbf{B}(t,\ee^{-\ii\pi/4}R)=
    2\sqrt{\pi t\epsilon}\ee^{-\pi P/4}R^{-\ii\alpha/\epsilon}\left(1,O(R^{-2}),\dots,O(R^{-2})\right),\quad R\to+\infty,
\end{equation}
where $P$ is given by \eqref{eq:P-def}. To analyze the remaining rows of $\mathbf{B}(t,\ee^{-\ii\pi/4}R)$ in the same limit, we first observe that by invertible row operations and redefinition of the branch cuts $\Gamma_j$, $j=1,\dots,N$ to diagonal rays $\Gamma'_j:=p_j+\ee^{3\pi\ii/4}\mathbb{R}_+$, at the cost of a nonzero constant factor in $\det(\mathbf{B}(t,x))$ we may assume that the contour $C_{j-1}$ in row $j$ is replaced with a ``U-shaped'' path surrounding with positive orientation only the new branch cut $\Gamma'_{j-1}$ (unless the index $j-1$ is exceptional, in which case we can simply take $C_{j-1}$ itself to be the ray $p_{j-1}+\ee^{3\pi\ii/4}\mathbb{R}_+$ oriented away from $p_{j-1}$). We observe that all the dependence on $x=R\ee^{-\ii\pi/4}$ in $\ee^{-\ii h(z)/\epsilon}$ enters via the exponential factor 
\begin{equation}
    \exp\left(-\ii \frac{(z-R\ee^{-\ii\pi/4})^2}{4t\epsilon}\right)= 
    \exp\left(-\frac{R^2}{4t\epsilon}\right)\exp\left(R\frac{\ee^{\ii\pi/4}z}{2t\epsilon}\right)\exp\left(-\frac{\ii z^2}{4t\epsilon}\right)
\end{equation}
and aside from the first factor on the right-hand side that does not depend on $z$, the only dependence on $R>0$ in $\ee^{-\ii h(z)/\epsilon}$ arises from the second factor.  This factor (having an exponent that is proportional to $R$) decays exponentially in the direction $\arg(z)=3\pi/4$ of the modified branch cuts and it is maximized at the corresponding branch point.  The relevant steepest-descent procedure for such situations is based on Watson's Lemma combined with repeated integration by parts to obtain integrability at the branch point (see \cite[Section 4.8]{AAA}). Using this method, one finds that row $j=2,3,\dots,N+1$ of $\mathbf{B}(t,\ee^{-\ii\pi/4}R)$ can be written in the form
\begin{multline}
    \mathbf{e}_{j}^\top\mathbf{B}(t,\ee^{-\ii\pi/4}R) = C_j(t,\epsilon)\exp\left(-\frac{R^2}{4t\epsilon}\right)\exp\left(R\frac{\ee^{\ii\pi/4}p_{j-1}}{2t\epsilon}\right)R^{-\ii c_{j-1}/\epsilon}\\
    \cdot\left(O(R^{-1}),\dots,O(R^{-1}),1,O(R^{-1}),\dots,O(R^{-1})\right),\quad R\to+\infty,
\end{multline}
where $C_j(t,\epsilon)\neq 0$ is a constant, and the largest entry is in the $j^\mathrm{th}$ column.

The determinant of $\mathbf{B}(t,x)$ is therefore dominated in this limit by the product of the diagonal entries, which is nonzero for sufficiently large $R>0$.
\end{proof}
With this result, we can now prove 
\begin{proposition}\label{lem:det-B-not-zero}
We have $\det(\mathbf{B}(t,x))\neq 0$ for every $x\in \mathbb{R}$ and $t,\epsilon>0$.
\end{proposition}

\begin{proof}
Given $t,\epsilon>0$, it follows from Lemma~\ref{lem:isolated-zeros}  that $x\mapsto\det(\mathbf{B}(t,x))$ is an entire function having only isolated zeros. 
Given the form of the related determinant $\det(\overline{\mathbf{B}}(t,x))$ in~\eqref{eq:Btilde}, we deduce that for $t,\epsilon>0$ fixed, $x\mapsto\det(\overline{\mathbf{B}}(t,x))$ is  analytic on $\mathbb{R}$ with only isolated zeros, exactly the same zeros as for $x\mapsto\det(\mathbf{B}(t,x))$. Suppose that $x_0\in\mathbb{R}$ is such that $\det(\mathbf{B}(t,x_0))=0$.  Then, $x=x_0$ is also an isolated zero of the analytic function $x\mapsto\det(\overline{\mathbf{B}}(t,x))$, so there exists an analytic function $x\mapsto\mu(x)$ depending parametrically on $t$ with $\mu(x_0)\neq 0$ and $k\in\mathbb{N}^*$ such that for $x$ in the vicinity of $x_0$, 
\begin{equation}
\det(\overline{\mathbf{B}}(t,x))=\mu(x) (x-x_0)^k.
\end{equation}
Therefore, given~\eqref{eq:B-to-Btilde}, this means that for $x\neq x_0$ in the vicinity of $x_0$, there holds
\begin{equation}
\Pi u(t,x)
    =\ii\epsilon\partial_x\log(\det(\overline{\mathbf{B}}(t,x)))=\frac{\ii\epsilon k}{x-x_0}+\ii\epsilon\frac{\partial_x\mu(x)}{\mu(x)}.
\end{equation}
The right-hand side is clearly not in $L^2_\mathrm{loc}(\mathbb{R})$ in the vicinity of $x_0$. However, global well-posedness of the Cauchy problem associated to~\eqref{eq:BO} in $H^2(\mathbb{R})$ was proven in~\cite{Saut79}, ensuring that $x\mapsto\Pi u(t,x)$ is a well-defined function in $H^2(\mathbb{R})$.  This is therefore a contradiction with the assumption that $\det(\mathbf{B}(t,x_0))=0$.
\end{proof}

This completes the proof of Theorem~\ref{thm:inversion-formula}.

\section{Long time asymptotics for specific initial data}\label{section:weak-asymptotics}
The next two sections of the paper are concerned with the proof of Theorem~\ref{thm:minussoliton}.
In this section, we apply the exact solution formula given in Theorem~\ref{thm:inversion-formula} in the special case that $\epsilon=1$ and the initial data is given by \eqref{eq:minus-soliton}, studying the solution in the limit $t\to+\infty$.  First we establish convergence of a renormalized version of $\Pi u(t,x)$ in frames of reference $x=ct$ with constant negative velocities $c=2y<0$.

\begin{proposition}[Locally uniform asymptotics]\label{prop:weak-asymptotics}
Let $u(t,x)$ denote the solution of \eqref{eq:BO} with $\epsilon=1$ and initial data $u_0(x)=-2/(1+x^2)$.  Let $\phi(t,y)$ be defined for $t>0$ and $y<0$ by
\begin{equation}
    \phi(t,y):=\ee^{-\ii\pi/4}\sqrt{4\pi t}\ee^{\ii t y^2}\Pi u(t,2ty).
\label{eq:renormalized-Pi-u}
\end{equation}
Then $\phi(t,y)\to\widehat{\psi}(-y)$ as $t\to+\infty$ locally uniformly for $y<0$,
where
\begin{equation}
\label{eq:widehatpsi}
\widehat{\psi}(\lambda)
	:=\frac{2\pi \ee^{\lambda}}{ \mathrm{Ei}(2\lambda)+\ii\pi},\quad\lambda>0,
\end{equation}
and $\mathrm{Ei}(\diamond)$ is given by \eqref{eq:Ei}.
\end{proposition}

\begin{proof}
The initial data $u_0(x)=-2/(1+x^2)$ is of the form \eqref{eq:rationalIC} with $N=1$, $p_1=\ii$, and $c_1=\ii$.  Since $\ii c_1/\epsilon=-1$, the (only) index $1$ is exceptional in the sense of Definition~\ref{def:C}.  Therefore, in this situation, $\ee^{-\ii h(z)}$ is single-valued and given explicitly by
\begin{equation}
\ee^{-\ii h(z)}=\ee^{-\ii(z-x)^2/(4t)}\frac{z-\ii}{z+\ii}.
\label{eq:exponential}
\end{equation}
Thanks to Theorem~\ref{thm:inversion-formula}, we then have the following explicit representation of the projection $\Pi u(t,x)$:
\begin{equation}
\Pi u(t,x)=\frac{N(t,x)}{D(t,x)},
\end{equation}
\begin{equation}
N(t,x)
	=\begin{vmatrix}
\displaystyle \int_{C_{0}}\ee^{-\ii h(z)}u_0(z)\,\dd z & 	\displaystyle \int_{C_{0}}\ee^{-\ii h(z)}\frac{\dd z}{z-\ii} \\
\displaystyle \int_{C_1}\ee^{-\ii h(z)}u_0(z)\,\dd z & 	\displaystyle \int_{C_1}\ee^{-\ii h(z)}\frac{\dd z}{z-\ii} \\
\end{vmatrix},
\end{equation}
\begin{equation}
D(t,x)
	=\begin{vmatrix}
\displaystyle \int_{C_{0}}\ee^{-\ii h(z)}\,\dd z & \displaystyle	\int_{C_{0}}\ee^{-\ii h(z)}\frac{\dd z}{z-\ii} \\
\displaystyle \int_{C_1}\ee^{-\ii h(z)}\,\dd z & \displaystyle	\int_{C_1}\ee^{-\ii h(z)}\frac{\dd z}{z-\ii} \\
\end{vmatrix}.
\end{equation}
Using \eqref{eq:exponential} and noting that 
\begin{equation}
    u_0(z)\frac{z-\ii}{z+\ii}=\frac{-2}{(z+\ii)^2},
\end{equation}
we may write the determinant $N(t,x)$ in the form
\begin{equation}\label{def:N}
N(t,x)
=-2\begin{vmatrix}J_0(t,x) & I_0(t,x)\\J_1(t,x) & I_1(t,x)\end{vmatrix},
\end{equation}
where for $j=0,1$,
\begin{equation}
    I_j(t,x):=\int_{C_j}\ee^{-\ii (z-x)^2/(4t)}\frac{\dd z}{z+\ii},\quad J_j(t,x):=\int_{C_j}\ee^{-\ii(z-x)^2/(4t)}\frac{\dd z}{(z+\ii)^2}.
\end{equation}
Similarly, using \eqref{eq:exponential} and adding $2\ii$ times the second column to the first, we write $D(t,x)$ in the form
\begin{equation}\label{def:D}
D(t,x)=\begin{vmatrix}K_0(t,x) & I_0(t,x)\\K_1(t,x) & I_1(t,x)\end{vmatrix},
\end{equation}
where for $j=0,1$,
\begin{equation}
    K_j(t,x):=\int_{C_j}\ee^{-\ii (z-x)^2/(4t)}\,\dd z.
\end{equation}
Recalling that 
\begin{itemize}
\item $C_0$ is a contour going from $z=\infty$ with $\arg(z)=3\pi/4$ to $z=\infty$ with $\arg(z)=-\pi/4$, passing above the pole at $z=-\ii$;
\item  $C_1$ is a contour in the upper half-plane originating from $z=\infty$ with $\arg(z)=3\pi/4$ and terminating at $z=\ii$ (as the index $1$ is exceptional),
\end{itemize}
we shall now calculate the integrals $I_0, I_1, J_0, J_1, K_0, K_1$.  To be concrete, we take $C_0$ to be the straight line through $z=\ii$ with $\arg(z-\ii)=3\pi/4$ or $-\pi/4$, and $C_1$ to be the ray $\arg(z-\ii)=3\pi/4$.

First, let us consider $I_0, J_0$ and $K_0$. The simplest is $K_0(t,x)$, which is an exact Gaussian integral with value independent of $x$:
\begin{equation}
    K_0(t,x)=\ee^{-\ii\pi/4}2\sqrt{\pi t}.
\end{equation}
Now let $y<0$ and set $x=2ty<0$. To calculate $I_0(t,2ty)$ and $J_0(t,2ty)$, we define the contour $C_0'$ parallel to $C_0$ but passing through the point $z=2ty$. When $t$ is large enough, $C_0'$ passes below the pole at $z=-\ii$ so the concatenation of the contour $C_0$ (following its orientation) and $C_0'$ (following the opposite direction) encircles in the clockwise sense the pole located at $z=-\ii$. Therefore, by residues,
\begin{equation}
    I_0(t,2ty)=-2\pi\ii \ee^{-\ii t y^2}\ee^y\ee^{\ii/(4t)}+\int_{C_{0}'}\ee^{-\ii(z-2ty)^2/(4t)}\frac{\dd z}{z+\ii},
\end{equation}
\begin{equation}
    J_0(t,2ty)=-2\pi\ii\left(\ii y-\frac{1}{2t}\right)\ee^{-\ii ty^2}\ee^y\ee^{\ii/(4t)}+\int_{C_{0}'}\ee^{-\ii(z-2ty)^2/(4t)}\frac{\dd z}{(z+\ii)^2}.
\end{equation}
The parameterization of $C_0'$ by $z=2ty+\ee^{-\ii\pi/4}s$ for $s\in \mathbb{R}$ leads for $p=1,2$ to
\begin{equation}
    \int_{C_0'}\ee^{-\ii (z-2ty)^2/(4t)}\frac{\dd z}{(z+\ii)^p} =\ee^{\ii\pi (p-1)/4}\int_\mathbb{R}\ee^{-s^2/(4t)}\frac{\dd s}{(s+\ee^{\ii\pi/4}(2ty+\ii))^p}.
\end{equation}
Given $y<0$ and $t>0$ sufficiently large we have
\begin{equation}
    \inf_{s\in\mathbb{R}} |s+\ee^{\ii\pi/4}(2ty+\ii)| =-\sqrt{2}ty-\frac{1}{\sqrt{2}}\ge -ty>0.
\end{equation}
Hence for both $p=1,2$,
\begin{equation}
    \left|\int_{C_0'}\ee^{-\ii (z-2ty)^2/(4t)}\frac{\dd z}{(z+\ii)^p}\right|\le \frac{1}{(-ty)^p}\int_\mathbb{R}\ee^{-s^2/(4t)}\,\dd s = \frac{2\sqrt{\pi t}}{(-ty)^p} = O(t^{1/2-p}),\quad t\to+\infty.
\end{equation}
Since $\ee^{\ii/(4t)}=1+O(t^{-1})$ as $t\to+\infty$, we conclude that in this limit,
\begin{equation}
\begin{split}
I_0(t,2ty)
	&=-2\pi\ii \ee^{-\ii ty^2+y}+O\left(\frac{1}{\sqrt{t}}\right),\\
J_0(t,2ty)
	&=2\pi y\ee^{-\ii ty^2+y}+O\left(\frac 1t\right).
    \end{split}
\end{equation}
The error terms in each case are locally uniform for $y<0$.

To calculate $I_1(t,2ty)$, $J_1(t,2ty)$, and $K_1(t,2ty)$, first observe that if $C_1'$ is the ray $\arg(z-\ii)=-\pi/4$ oriented toward $z=\ii$, then
\begin{equation}
    \begin{split}
        I_1(t,2ty)&=I_0(t,2ty) + \int_{C_1'}\ee^{-\ii (z-2ty)^2/(4t)}\frac{\dd z}{z+\ii}\\
        &=I_0(t,2ty) + \ee^{-\ii ty^2}\int_{C_1'}\ee^{-\ii z^2/(4t)}\ee^{\ii yz}\frac{\dd z}{z+\ii},\\
        J_1(t,2ty)&=J_0(t,2ty) + \int_{C_1'}\ee^{-\ii (z-2ty)^2/(4t)}\frac{\dd z}{(z+\ii)^2}\\
        &=J_0(t,2ty) +\ee^{-\ii ty^2}\int_{C_1'}\ee^{-\ii z^2/(4t)}\ee^{\ii yz}\frac{\dd z}{(z+\ii)^2},\\
        K_1(t,2ty)&=K_0(t,2ty) +\int_{C_1'}\ee^{-\ii (z-2ty)^2/(4t)}\,\dd z\\
        &=K_0(t,2ty)+\ee^{-\ii ty^2}\int_{C_1'}\ee^{-\ii z^2/(4t)}\ee^{\ii yz}\,\dd z.
    \end{split}
\end{equation}
The reason for expressing $I_1(t,2ty)$, $J_1(t,2ty)$, and $K_1(t,2ty)$ in this way is that in each case the integral over $C_1'$ on the second line is amenable to the Lebesgue dominated convergence theorem to allow taking the limit $t\to+\infty$ under the integral sign.  Hence
\begin{equation}
    \begin{split}
        I_1(t,2ty)&=I_0(t,2ty) + \ee^{-\ii ty^2}F(y)\int_{C_1'}\ee^{\ii y z}\frac{\dd z}{z+\ii} + o(1),\\
        J_1(t,2ty)&=J_0(t,2ty) + \ee^{-\ii t y^2}\int_{C_1'}\ee^{\ii y z}\frac{\dd z}{(z+\ii)^2} + o(1),\\
        K_1(t,2ty)&=K_0(t,2ty) + \ee^{-\ii t y^2}\int_{C_1'}\ee^{\ii yz}\,\dd z + o(1),
    \end{split}
\end{equation}
where the $o(1)$ error terms are all locally uniform for $y<0$, and where
\begin{equation}
F(y):=\int_{C_1'}\ee^{\ii y z}\frac{\dd z}{z+\ii},\quad y<0.
    \label{eq:F-of-y-def}
\end{equation}
Of course we have explicitly
\begin{equation}
    \int_{C_1'}\ee^{\ii y z}\,\dd z = \frac{\ee^{-y}}{\ii y},
\end{equation}
and integrating by parts,
\begin{equation}
    \int_{C_1'}\ee^{\ii y z}\frac{\dd z}{(z+\ii)^2} = -\frac{\ee^{-y}}{2\ii} +\ii yF(y).
\end{equation}

Using these results in \eqref{def:N} and \eqref{def:D} and subtracting the first row from the second in each case gives
\begin{equation}
\begin{split}
    N(t,2ty)&=-2\begin{vmatrix}J_0(t,2ty) & I_0(t,2ty)\\J_1(t,2ty)-J_0(t,2ty) & I_1(t,2ty)-I_0(t,2ty)\end{vmatrix}\\ 
    &=-2\begin{vmatrix}2\pi y\ee^{-\ii ty^2+y} + o(1) & -2\pi\ii\ee^{-\ii ty^2+y}+o(1)\\
    \frac{1}{2}\ii\ee^{-\ii ty^2-y}+\ii y\ee^{-\ii ty^2}F(y)+o(1) & \ee^{-\ii ty^2}F(y)+o(1)\end{vmatrix}\\
    &= 2\pi\ee^{-2\ii ty^2} + o(1)
\end{split}
\end{equation}
and
\begin{equation}
\begin{split}
    D(t,2ty)&=\begin{vmatrix}K_0(t,2ty) & I_0(t,2ty)\\K_1(t,2ty)-K_0(t,2ty) & I_1(t,2ty)-I_0(t,2ty)\end{vmatrix}\\ &= \begin{vmatrix}\ee^{-\ii\pi/4}2\sqrt{\pi t} & -2\pi\ii\ee^{-\ii t y^2+y} + o(1)\\
    -\ii y^{-1}\ee^{-\ii ty^2-y} + o(1) & \ee^{-\ii ty^2}F(y) + o(1)\end{vmatrix}\\
    &=2\sqrt{\pi t}\ee^{-\ii\pi/4}\ee^{-\ii ty^2}F(y) + o(t^{1/2}).
    \end{split}
\end{equation}
Therefore, we deduce that 
\begin{equation}
    \Pi u(t,2ty)=\frac{N(t,2ty)}{D(t,2ty)}=\sqrt{\frac{\pi}{t}}\frac{\ee^{\ii\pi/4}\ee^{-\ii ty^2}}{F(y)} + o(t^{-1/2})
    \label{eq:Pi-u-asymp}
\end{equation}
as $t\to+\infty$, locally uniformly for $y<0$.  

It only remains to calculate the integral $F(y)$,
which can be expressed in terms of the special function $\mathrm{Ei}(\diamond)$ as follows.  First, without changing the value of $F(y)$, we may rotate the contour $C_1'$ clockwise about the fixed endpoint $z=\ii$ so that it consists of, for arbitrarily small $\delta>0$, the concatenation of:
\begin{itemize}
    \item the negative imaginary axis from $z=-\ii\infty$ to $z=-\ii-\ii\delta$;
    \item the circular arc $z=-\ii + \delta\ee^{\ii\theta}$ with $\theta$ increasing from $-\pi/2$ to $\pi/2$; 
    \item the imaginary axis from $z=-\ii+\ii\delta$ to $z=\ii$.
\end{itemize}
Then taking the limit $\delta\downarrow 0$, the contribution from the circular arc yields a contribution of $\pi\ii\ee^{y}$, and with the real parametrization $z=\ii w$ the remainder is a principal-value integral.  Thus:
\begin{equation}
F(y)=\pi\ii\ee^y + \text{P.V.}\int_{-\infty}^1\ee^{-yw}\frac{\dd w}{w+1}.
\end{equation}
Finally, making the substitution $w=-1-s/y$,
\begin{equation}
    F(y)=\ee^y\left(\pi\ii +\text{P.V.}\int_{-\infty}^{-2y}\frac{\ee^s}{s}\,\dd s\right) = \ee^y\left(\pi\ii+\mathrm{Ei}(-2y)\right),\quad y<0,
    \label{eq:F-Ei}
\end{equation}
where we used the definition \eqref{eq:Ei}.  Using this in \eqref{eq:Pi-u-asymp} and
referring to \eqref{eq:renormalized-Pi-u} and \eqref{eq:widehatpsi}, the proof is complete.
\end{proof}

Note that according to \eqref{eq:widehatpsi} and \cite[Eqn.\@ 6.12.2]{DLMF},
\begin{equation}
    \widehat{\psi}(\lambda)=4\pi\lambda\ee^{-\lambda}(1+o(1)),\quad\lambda\to+\infty
\end{equation}
and according to \cite[Eqn.\@ 6.6.1]{DLMF},
\begin{equation}
    \widehat{\psi}(\lambda)=\frac{2\pi}{\ln(\lambda)} + O(\lambda/\ln(\lambda)),\quad\lambda\downarrow 0.
\end{equation}
As $\mathrm{Ei}(2\lambda)$ is real for $\lambda>0$, one has $\widehat{\psi}\in L^2(0,\infty)$ and hence $\|\widehat{\psi}(-\diamond)\|_{L^2(-\infty,0)}<\infty$.  Also, using conservation of the $L^2$ norm under \eqref{eq:BO} a scaling argument shows that $\|\phi(t,\diamond)\|_{L^2(-\infty,0)}$ is bounded uniformly in $t$.  With these facts we can then prove the following result.

\begin{corollary}[Weak convergence in $L^2(-\infty,0)$]
As $t\to+\infty$, we have the weak convergence in $L^2(-\infty,0)$:
\begin{equation}
    \phi(t,\diamond)\rightharpoonup\widehat{\psi}(-\diamond).
\end{equation}
\label{cor:weak-convergence}
\end{corollary}
\begin{proof}
First, using that the convergence of $\phi(t,\diamond)$ to $\widehat{\psi}(-\diamond)$ is locally uniform, we deduce that it is uniform on every compact subset of $(-\infty,0)$. Hence for each $v\in \mathcal{C}_\mathrm{c}^\infty(-\infty,0)$, we can see that the convergence is uniform on the support of $v$ so that
\begin{equation}
    |\langle v,\phi(t,\diamond)-\widehat{\psi}(-\diamond)\rangle|\to 0,\quad t\to\infty.
    \label{eq:smooth-weak}
\end{equation}
Then, using that $\phi(t,\diamond)$ and $\widehat{\psi}(-\diamond)$ are bounded in $L^2(-\infty,0)$ uniformly as $t\to+\infty$, there is a constant $M>0$ such that for each given $g\in L^2(-\infty,0)$ and $v\in \mathcal{C}^\infty_\mathrm{c}(-\infty,0)$,
\begin{equation}
     |\langle g,\phi(t,\diamond)-\widehat{\psi}(-\diamond)\rangle|
     \leq  |\langle v,\phi(t,\diamond)-\widehat{\psi}(-\diamond)\rangle|+ M\|v-g\|_{L^{2}(-\infty,0)},\quad t>0.
\end{equation}
According to \eqref{eq:smooth-weak}, the $\limsup$ as $t\to+\infty$ of the right-hand side is bounded by $M\|g-v\|_{L^2(-\infty,0)}$, which can be made arbitrarily small for given $g\in L^2(-\infty,0)$ by choice of $v\in \mathcal{C}_\mathrm{c}^\infty(-\infty,0)$ thanks to a density argument. This is enough to conclude that
\begin{equation}
     |\langle g,\phi(t,\diamond)-\widehat{\psi}(-\diamond)\rangle|\to 0,\quad t\to +\infty
\end{equation}
holds for every $g\in L^2(-\infty,0)$.
\end{proof}

\section{Spectral theory of the Lax operator}\label{section:spectral-Lax}

In this section, we develop a spectral theory for the Lax operator $L_{u_0}:=-\ii\epsilon\partial_x-T_{u_0}$ if $\langle x\rangle u_0\in L^2(\R )$ (recall, from \eqref{eq:Toeplitz}, that $T_{u_0}f=\Pi(u_0f)$). In this case, it is well known --- see e.g. \cite[Prop. 2.4]{Sun2020} --- that the essential spectrum of $L_{u_0}$ is $[0,+\infty )$ and that the eigenvalues of $L_{u_0}$ are strictly negative and simple. Consequently, if $u_0\leq 0$, then $L_{u_0}$ is a positive operator  and hence has no eigenvalues. We denote by $L^p_+(\R )$ the space of $L^p $ functions on $\R $ with Fourier transform supported in $[0,\infty )$. 

\begin{remark}
In this section, we develop the general theory assuming that $\epsilon=1$; however note that if $u(t,x)$ is the solution of the Benjamin-Ono equation in the general form \eqref{eq:BO} with $u(t,0)=u_0(x)$, then $v(t,x):=\epsilon^{-1}u(\epsilon^{-1}t,x)$ is the solution of the same equation with $\epsilon=1$ and with initial data $v(0,x)=v_0(x)=\epsilon^{-1}u_0(x)$.  Similarly, if we use the notation $L^{(\epsilon)}_{u_0}:=-\ii\epsilon\partial_x-T_{u_0}$ to indicate the dependence on $\epsilon$, then $L^{(\epsilon)}_{u_0}=\epsilon L^{(1)}_{v_0}$.  Therefore it is sufficient to assume here that $\epsilon=1$. Also, to keep the notation simple, we write $u$ in place of $u_0$ throughout this section.  
\end{remark}

\subsection{A distorted Fourier transform}
\label{sec:distortedFT}
Our immediate main goal is the following theorem.
\begin{theorem}[Distorted Plancherel theorem]\label{spectralLu}
Denote by $m_-(\cdot,\lambda )$ the unique solution in $L^\infty_+(\R )$ of the equation
\begin{equation}
    (L_u-\lambda)m=0
\end{equation}
such that $\ee^{-\ii\lambda x}m_-(x,\lambda )\td_x,{-\infty}1$,
and, for every $f\in L^1_+\cap L^2_+$, 
\begin{equation}
    \widetilde f(\lambda ):=\int_\R f(x)m_-(x,\lambda)^*\,\dd x.
\label{eq:distortedFT}
\end{equation}
Then the map $f\mapsto  \widetilde f$ extends to a bounded map from $L^2_+(\mathbb{R})$ to $L^2(0,\infty )$, and, for every $f\in L^2_+(\mathbb{R})$, the spectral measure $\mu_f$ of $f$ for $L_u$ is given by
\begin{equation}
    \dd\mu_f(\lambda )=\sum_{j}|\langle f\vert \phe_j\rangle |^2\delta_{\lambda _j}+\frac{1}{2\pi}\vert \widetilde f(\lambda )\vert ^2\, \dd\lambda,
\end{equation}
where the $\lambda_j$ are the eigenvalues of $L_u$, and the $\phe_j$ are corresponding normalized eigenfunctions.
In particular,
\begin{equation}
    \Vert f\Vert_{L^2}^2=\sum_{j}|\langle f\vert \phe_j\rangle |^2+\frac{1}{2\pi}\int_0^{+\infty} \vert \widetilde f(\lambda )\vert ^2\, \dd\lambda.
\end{equation}
In the special case $f=\Pi u$, we infer
\begin{equation}
    \Vert \Pi u\Vert_{L^2}^2=\sum_j 2\pi \vert \lambda_j\vert  +\frac{1}{2\pi}\int_0^{+\infty}\vert \widetilde {\Pi u}(\lambda)\vert ^2\, \dd\lambda.
\end{equation}
\end{theorem}
\begin{corollary}\label{spectralLuplus}
If moreover $L_u$ is a positive operator --- in particular if $u\leq 0$ --- then 
\begin{equation}
    \Vert \Pi u\Vert_{L^2}^2=\frac{1}{2\pi}\int_0^{+\infty}\vert \widetilde {\Pi u}(\lambda)\vert ^2\, \dd\lambda.
\end{equation}
\end{corollary}
\begin{proof}
We shall make extensive use of the following lemma.
\begin{lemma}\label{xT}
Let $u(\diamond)\in L^2(\R )$ be real valued, such that also $\diamond u(\diamond)\in L^2(\R)$.  
Then, for every $m\in L^\infty_+(\R )$, 
\begin{equation}
    T_um(x)=g (x) +\frac{c} {x-\ii}
\end{equation}
with $g \in L^1(\R )$ and $c\in \C$. More precisely,
\begin{equation}
    g(x) :=\frac{T_{(\diamond-\ii)u}m(x)}{x-\ii}\ ,\ c:=\frac{\ii}{2\pi}\int_\R u(x)m(x)\, \dd x.
\end{equation}
\end{lemma}
\begin{proof}
This is a consequence of the following general fact. For every $f(\diamond)\in L^2(\R )$ such that $\diamond f(\diamond)\in L^2(\R)$.
\begin{equation}
    x\Pi f(x)=\Pi (\diamond f(\diamond))(x)+\frac{\ii}{2\pi}\int_\R f(x)\, \dd x,
\end{equation}
which is obvious via the Fourier transformation. Just apply this identity to $f=um$.
\end{proof}
We now come to the first result.
\begin{proposition}\label{geneigen}
Given $\lambda >0$, the vector space 
$$\mathscr M(u,\lambda ):=\{ m\in L^\infty_+(\R ): (L_u-\lambda )m=0\}$$
is one dimensional, and the linear forms 
$$\ell_{\pm }(m,\lambda )=\lim_{x\to \pm \infty} m(x){\rm e}^{-\ii\lambda x}$$
are well defined and not identically zero on $\mathscr M(u,\lambda )$. 
\end{proposition}
\begin{proof}
Let $m\in \mathscr M(u,\lambda )$. Then
\begin{equation}
    \frac{\dd}{\dd x}(\ee^{-\ii\lambda x}m(x))=\ii\ee^{-\ii\lambda x}T_um(x)=\ii\ee^{-\ii\lambda x}g(x)+\ee^{-\ii\lambda x}\frac{\ii c}{x-\ii},
\end{equation}
where we have used Lemma \ref{xT}. Since $g\in L^1(\R )$, so is $x\mapsto \ee^{-\ii\lambda x}g(x)$.  Moreover, since $\lambda \ne 0$, via integration by parts, one can show that
the generalized integrals
\begin{equation}
    \int_{-\infty} ^0\frac{\ee^{-\ii\lambda x}}{x-\ii}\, \dd x\ ,\ \int_0^{+\infty} \frac{\ee^{-\ii\lambda x}}{x-\ii}\, \dd x
\end{equation}
are convergent. Hence $\ell_+ (m,\lambda )$ and $\ell_-(m,\lambda)$ are well defined.\\

Let us now prove that the kernel of $\ell_\pm (\cdot,\lambda)$ on $\mathscr M(u,\lambda )$ is $\{ 0\}$.
First of all, we observe that
\begin{multline}
    \frac{\dd}{\dd x}|m(x)|^2=2{\rm Re}( m(x)^*m'(x))=-2{\rm Im}( m(x)^*T_um(x)) \\
    =-2{\rm Im}( m(x)^*g(x))-2{\rm Im}(m(x)^*c(x-\ii)^{-1})\ .
\end{multline}
The following identity holds for the integral of the first term in the right hand side,
\begin{equation}
    \int_\R  m(x)^*g(x)\, \dd x=\left \langle T_{(\diamond-\ii)u}m, \frac{m}{\diamond+\ii}\right \rangle=\left \langle (\diamond-\ii)um, \frac{m}{\diamond+\ii}\right \rangle=\int_\R u(x)|m(x)|^2\, \dd x \in \R,
\end{equation}
hence its imaginary part vanishes. On the other hand,
\begin{eqnarray*}
\int_{-R}^R \frac{ m(x)^*}{x-\ii}\, \dd x&=&\int_{-R}^R \frac{m(x)^*\ee^{\ii\lambda x}}{x-\ii}\ee^{-\ii\lambda x}\, \dd x\\
&=&O(R^{-1})+(\ii\lambda )^{-1}\int_{-R}^R \frac{\dd}{\dd x}\left [ \frac{ m(x)^*\ee^{\ii\lambda x}}{x-\ii}\right ] \ee^{-\ii\lambda x}\, \dd x\\
&=&O(R^{-1})-\lambda ^{-1}\int_{-R}^R  \frac{T_um(x)^*}{x-\ii} \, \dd x -(\ii\lambda)^{-1}\int_{-R}^R \frac{ m(x)^*}{(x-\ii)^2}\, \dd x\\
&=&O(R^{-1})-\lambda^{-1}\langle (\diamond-\ii)^{-1}, T_um\rangle -(\ii\lambda)^{-1}\langle (\diamond-\ii)^{-2}, m\rangle \\
&=&O(R^{-1})\ ,
\end{eqnarray*}
since $(x-\ii)^{-1}\in L^2_-(\R)\perp T_um\in L^2_+(\R )$, and $(x-\ii)^{-2}\in L^1_-(\R )\perp m\in L^\infty_+(\R )$.

We conclude that $|m(R)|^2-|m(-R)|^2\to 0$ as $R\to +\infty $, hence $|\ell_+(m,\lambda )|^2=|\ell_-(m,\lambda )|^2$, and the kernels of these two linear forms on $\mathscr M(u,\lambda )$ coincide. 
Let $m$ be in this kernel. Then, using Lemma \ref{xT},
\begin{equation}
    m(x)=\ii\int_{-\infty}^x \ee^{\ii\lambda (x-y)}T_um(y)\, \dd y=\ii\int_{-\infty}^x \ee^{\ii\lambda (x-y)}g(y)\, \dd y+\ii c\int_{-\infty}^x \ee^{\ii\lambda (x-y)}(y-\ii)^{-1}\, \dd y.
\end{equation}
The first term on the right-hand side is in $L^2(-\infty, 0)$ because of the Hardy inequality and $|g(y)|\leq (1+|y|)^{-1}h(y)$ with $h\in L^2(\R )$. As for the second term, an integration by parts implies that it is $O((1+|x|)^{-1})$. Hence we conclude that $m\in L^2(-\infty, 0)$. Similarly, we have $m\in L^2(0,+\infty )$. We infer that $m\in L^2_+(\R )$ and is therefore an eigenvector of $L_u$ with the positive eigenvalue $\lambda $. This implies $m=0$.

In order to complete the proof, we just need to establish that the vector space $\mathscr M(u,\lambda )$ is not trivial.
Consider the operator $K_{u,\lambda}:L^\infty_+\to L^\infty _+$ defined by
\begin{equation}
    K_{u,\lambda}m(x)=\int_{-\infty}^x \ee^{\ii\lambda (x-y)}T_um(y)\, \dd y=\int_0^{+\infty} \ee^{\ii\lambda t}T_um(x-t)\, \dd t.
\end{equation}
Because of Lemma \ref{xT}, we have
\begin{equation}
    K_{u,\lambda}m(x)=\int_{-\infty}^x \ee^{\ii\lambda(x-y)}g(y)\, \dd y+c\int_{-\infty}^x\ee^{\ii\lambda(x-y)}\frac{\dd y}{y-\ii}.
\end{equation}
The first integral is $L^\infty $ because $g\in L^1$, while the second integral is shown to be $L^\infty $ --- with finite limits at $\pm \infty $ --- after an integration by parts. We claim 
\begin{lemma}\label{Kcomp}
    $K_{u,\lambda}:L^\infty_+\to L^\infty _+$ is a compact operator.
\end{lemma}
\begin{proof}
It is enough to prove that $\overline K_{u,\lambda}$ defined by
\begin{equation}
    \overline K_{u,\lambda}m(x):=\int_{-\infty}^x \ee^{-\ii\lambda y}T_um(y)\, \dd y
\end{equation}
is compact. On the one hand, we have, if $\| m\|_{L^\infty}\leq 1$,
\begin{equation}
    |\overline K_{u,\lambda}m(x)-\overline K_{u,\lambda}m(x')|\leq |x-x'|^{1/2}\| u\|_{L^2}.
\end{equation}
On the other hand, using again Lemma \ref{xT}, 
\begin{align}
    |\overline K_{u,\lambda}m(x)|&=O(|x|^{-1/2}) ,\quad x\to -\infty,\\
    |\overline K_{u,\lambda}m(x)-\overline K_{u,\lambda}m(+\infty )|&=O(|x|^{-1/2}),\quad x\to +\infty.
\end{align}
Hence the image of the unit ball of $L^\infty _+$ by $\overline K_{u,\lambda}$ is uniformly small at $-\infty$, uniformly close to (bounded) constants at $+\infty $, and uniformly equicontinuous. It is therefore relatively compact in $L^\infty$.
\end{proof}

Having proved the compactness of $K_{u,\lambda}$, we notice that the kernel of $1-\ii K_{u,\lambda}$ precisely consists of functions in $\mathscr M(u,\lambda )$ such that 
\begin{equation}
    \ell_-(m,\lambda)=0,
\end{equation}
which just have been shown to be identically $0$. 
By the Riesz--Fredholm theorem, we conclude that $1-\ii K_{u,\lambda}$ is bijective from $L^\infty_+$ onto $L^\infty _+$. Now the function $m_-$ defined by
\begin{equation}
    m_-(x,\lambda )=(1-\ii K_{u,\lambda})^{-1}(\ee^{\ii\lambda \diamond})(x),
\end{equation}
belongs to $\mathscr M(u,\lambda )$ and satisfies $\ell_-(m_-,\lambda )=1$. The proof of Proposition \ref{geneigen} is complete.
\end{proof}

We now come to the description of the spectral measure $\mu_f$ of $f$. By the limiting absorption principle (or by functional calculus for the self--adjoint operator $L_u$), we know that, for the weak-$*$ topology of measures,
\begin{equation}
    \dd\mu_f(\lambda )=\frac{1}{2\pi\ii}\lim_{\e \to 0^+} \langle [(L_u-\lambda -\ii\e)^{-1}-(L_u-\lambda +\ii\e)^{-1}]f, f\rangle.
\end{equation}
Furthermore, we know that operator $L_u$ has purely point spectrum on $(-\infty ,0)$ and has no eigenvalue at $0$. Therefore, we just have to calculate the above weak-$*$ limit on $(0,+\infty )$. Furthermore, since $L^1_+\cap L^2_+$ is dense in $L^2_+$, we may assume that $f\in L^1_+\cap L^2_+$. We set
\begin{equation}
    g_\e^\pm:=(L_u-\lambda \mp\ii\epsilon)^{-1}f.
\end{equation}
Then with the standard notation $D_x=-\ii\partial_x$ (with Fourier multiplier $\xi$)
\begin{equation}
    g_\e^{\pm}=(1-(D_x-\lambda\mp\ii\e )^{-1}T_u)^{-1}(D_x-\lambda\mp\ii\e )^{-1}f
\end{equation}
with
\begin{align}
    (D_x-\lambda -\ii\e )^{-1}f(x)&=\ii\int_{-\infty}^x \ee^{(\ii\lambda -\e)(x-y)}f(y)\, \dd y, \\ 
    (D_x-\lambda +\ii\e )^{-1}f(x)&=-\ii\int_{x}^{+\infty} \ee^{(\ii\lambda +\e)(x-y)}f(y)\, \dd y.
\end{align}
In particular, since $f\in L^1_+$, $(D_x-\lambda \mp \ii\e )^{-1}f$ converge in $L^\infty _+$ to some limits that we denote, respectively, by $(D_x-\lambda \mp \ii 0 )^{-1}f$.
Furthermore,
\begin{equation}
    (D_x-\lambda -\ii\e )^{-1}T_u=\ii K_{u,\lambda}^{\epsilon}, \quad K_{u,\lambda }^{\epsilon} m(x):=\int_{-\infty}^x \ee^{(\ii\lambda -\epsilon)(x-y)}T_um(y)\, \dd y.
\end{equation}
By the same arguments as in the proof of Lemma \ref{Kcomp}, one proves that the family $K_{u,\lambda}^\e$ is uniformly compact as $\e \to 0^+$. Similarly, the family
$(D_x-\lambda +\ii\epsilon)^{-1}T_u:L^\infty_+\to L^\infty _+$ is uniformly compact as $\e \to 0^+$.

Now we claim that $g_\e ^\pm $ is uniformly bounded in $L^\infty$. Let us prove it for $g_\e ^+$, say. Otherwise, up to extracting a subsequence, we would have $\Vert g_\e ^+\Vert_{L^\infty }\to +\infty $. Setting 
\begin{equation}
    m_\e =\frac{g_\e ^+}{\Vert g_\e ^+\Vert_{L^\infty}},
\end{equation}
we would have
\begin{equation}
    m_\e =\frac{(D_x-\lambda-\ii\e)^{-1}f}{\Vert g_\e ^+\Vert_{L^\infty}}+\ii K_{u,\lambda}^\e m_\e.
\end{equation}
By the uniform compactness of the family $K_{u,\lambda }^\e$ as $\e \to 0^+$, we would infer that  --- after possibly extracting another subsequence --- $m_\e $ converges strongly in $L^\infty _+$ to some $m$ with
\begin{equation}
    m=\ii K_{u,\lambda }m,
\end{equation}
which precisely means that $m\in \ker \ell_-(\cdot,\lambda )$ and therefore,  in view of Proposition \ref{geneigen},  implies $m=0$. Since $\Vert m_\e \Vert_{L^\infty}=1$, this would contradict the strong convergence of $m_\e$ to $m$ in $L^\infty $.

We conclude that $g^\e _+$ is bounded in $L^\infty _+$, therefore it converges strongly in $L^\infty _+$ to  $g^+\in L^\infty _+$ given by
\begin{equation}
    g^+=(1-(D_x-\lambda-\ii0)^{-1}T_u)^{-1}(D_x-\lambda-\ii0)^{-1}f.
\end{equation}
Similarly, $g_\e ^-$ converges strongly to $g^-$ given by
\begin{equation}
    g^-=(1-(D_x-\lambda+\ii0)^{-1}T_u)^{-1}(D_x-\lambda+\ii0)^{-1}f.
\end{equation}
Notice that $m:=g^+-g^-$ belongs to $\mathscr M(u,\lambda )$. By Proposition \ref{geneigen},
\begin{equation}
    m(x)=\ell_-(m,\lambda )m_-(x,\lambda ).
\end{equation}
It remains to calculate 
\begin{equation}
    \ell_-(m,\lambda )=\lim_{x\to -\infty}{\rm e}^{-\ii\lambda x}(g^+(x)-g^-(x))=-\lim_{x\to -\infty}{\rm e}^{-\ii\lambda x}g^-(x).
\end{equation}
Multiplying by $ m_-(x,\lambda )^*$ the equation
\begin{equation}
    (L_u-\lambda )g^-=f,
\end{equation}
integrating from $-R$ to $R$ and recalling that $g^-(x)\to 0$ as $x\to +\infty $, we conclude
\begin{equation}
    \ii\lim_{x\to -\infty}{\rm e}^{-\ii\lambda x}g^-(x)=\int_{-\infty}^{+\infty}f(x) m_-(x,\lambda)^*\, \dd x.
\end{equation}
In other words,
\begin{equation}
    \ell_-(m,\lambda )=\ii\int_{-\infty}^{+\infty}f(x) m_-(x,\lambda)^*\, \dd x\ =: \ii\tilde f(\lambda )
\end{equation}
Consequently,  we get, on $(0,+\infty )$,
\begin{equation}
    \dd\mu_f(\lambda )=\frac{1}{2\pi}|\widetilde f(\lambda )|^2\, \dd\lambda.
\end{equation}
On the other hand, in view of formula (2.3) in \cite[Prop.\@ 2.4]{Sun2020}, we know that for each eigenvalue $\lambda_j<0$ with normalized eigenfunction $\varphi_j$,
\begin{equation}
    |\langle \Pi u, \varphi_j\rangle |^2=2\pi |\lambda_j|.
\end{equation}
This completes the proof of Theorem~\ref{spectralLu}.
\end{proof}

\subsection{Application to long-time asymptotics}
We now come back to the proof of Theorem \ref{thm:minussoliton} and consider $u_0$ given by \eqref{eq:minus-soliton}. Let us describe $m_-(x,\lambda)$ and $\alpha(\lambda):= \widetilde{\Pi u_0}(\lambda )$ in this case when the dispersion parameter is $\epsilon=1$. First we compute $m_-(x,\lambda)$.  The equation $(L_{u_0}-\lambda)m_-(x,\lambda)=0$ reads
\begin{equation}\frac{\dd m_-}{\dd x}(x,\lambda)+\frac{2\ii}{1+x^2}m_-(x,\lambda)-\ii\lambda m_-(x,\lambda)= \frac{c}{x-\ii}\ ,\end{equation}
where the constant $c$ has to be chosen so that $m_-(x,\lambda )$ extends holomorphically to ${\rm Im}(x)>0$. Equivalently,
\begin{equation}\frac \dd{\dd x} \left ({\rm e}^{-\ii\lambda x}\frac{x-\ii}{x+\ii}m_-(x,\lambda)\right ) = \frac{c}{x+\ii}{\rm e}^{-\ii\lambda x}\ .\end{equation}
Imposing the required boundary condition  
\begin{equation}{\rm e}^{-\ii\lambda x}m_-(x,\lambda)\td_x,{-\infty} 1\ ,\end{equation}
we obtain
\begin{equation}\frac{x-\ii}{x+\ii}{\rm e}^{-\ii\lambda x}m_-(x,\lambda)=1+c\int_{-\infty}^x \frac{{\rm e}^{-\ii\lambda z}}{z+\ii}\, \dd z,\end{equation}
where on the right-hand side we have a convergent improper integral.
The analyticity of $m(x,\lambda )$ near $x=\ii$ implies that the left hand side vanishes at $x=\ii$, so we obtain
\begin{equation}c=-\frac{1}{\displaystyle{\int_{-\infty}^\ii \frac{{\rm e}^{-\ii\lambda z}}{z+\ii}\, \dd z} }\,\end{equation}
where the path of integration lies in the closed upper half $z$-plane, and finally
\begin{equation}m_-(x,\lambda )=\frac{x+\ii}{x-\ii}\, {\rm e}^{\ii\lambda x}\  \frac{\displaystyle{\int_{-\infty}^\ii \frac{{\rm e}^{-\ii\lambda z}}{z+\ii}\, \dd z-\int_{-\infty}^x \frac{{\rm e}^{-\ii\lambda z}}{z+\ii}\, \dd z} }{\displaystyle{\int_{-\infty}^\ii \frac{{\rm e}^{-\ii\lambda z}}{z+\ii}\, \dd z} } =\frac{x+\ii}{x-\ii}\ee^{\ii\lambda x}\frac{\displaystyle\int_x^\ii\frac{\ee^{-\ii\lambda z}}{z+\ii}\,\dd z}{\displaystyle\int_{-\infty}^\ii\frac{\ee^{-\ii\lambda z}}{z+\ii}\,\dd z}   \ ,\end{equation}
where the integration path for the integral in the numerator can be taken to be a straight line segment from $z=x\in\mathbb{R}$ to $z=\ii$.  Since $\lambda>0$, the contour integral in the denominator can be related to $F(-\lambda)$ defined in \eqref{eq:F-of-y-def} by a residue calculation at $z=-\ii$:
\begin{equation}
\begin{split}
    \int_{-\infty}^\ii\frac{\ee^{-\ii\lambda z}}{z+\ii}\,\dd z &= -2\pi\ii\ee^{-\lambda} + F(-\lambda)\\
    &=\ee^{-\lambda}(-\ii\pi +\mathrm{Ei}(2\lambda)),\quad \lambda>0,
    \end{split}
\end{equation}
where in the second line we used \eqref{eq:F-Ei}.

Now we use the obtained expression for $m_-(x,\lambda)$ to calculate $\alpha(\lambda)$.  
Since $\Pi u_0(x)=-\ii/(x+\ii)$, the definition \eqref{eq:distortedFT} gives
\begin{equation}
\alpha(\lambda):=\widetilde{\Pi u_0}(\lambda)=\int_{-\infty}^{+\infty} \Pi u_0(x)m_-(x,\lambda)^*\, \dd x =\frac{n(\lambda){\rm e}^{\lambda}}{\mathrm{Ei}(2\lambda )+\ii\pi}\ ,\end{equation}
where 
\begin{equation}
    n(\lambda ):=-\ii\int_{-\infty}^{+\infty}\frac{x-\ii}{(x+\ii)^2}\ee^{-\ii\lambda x}\int_x^{-\ii}\frac{\ee^{\ii\lambda z}}{z-\ii}\,\dd z\,\dd x.
\end{equation}
Since $\lambda>0$, the path of integration in the improper outer integral can be replaced with a contour $C$ originating at $z=\infty$ with $\arg(z)=-3\pi/4$ and terminating at $z=\infty$ with $\arg(z)=-\pi/4$ that passes above the pole at $x=-\ii$.  
We then parametrize the path of integration in the inner integral by $z=x-r(x+\ii)$ for $0<r<1$ and exchange the integration order to get
\begin{equation}
\begin{split}
    n(\lambda)&=\ii\int_C\frac{x-\ii}{x+\ii}\int_0^1\frac{\ee^{-\ii\lambda r x}}{(1-r)x-\ii (1+r)}\ee^{\lambda r}\,\dd r\,\dd x \\ &= \ii
    \int_0^1\ee^{\lambda r}\int_{C}\frac{x-\ii}{(x+\ii)((1-r)x-\ii(1+r))}\ee^{-\ii\lambda r x}\,\dd x\,\dd r.
\end{split}
\end{equation}
The integral over $C$ can then be evaluated for each $r>0$ by a residue at $x=-\ii$:
\begin{equation}
\begin{split}
    \int_{C}\frac{x-\ii}{(x+\ii)((1-r)x-\ii(1+r))}\ee^{-\ii\lambda r x}\,\dd x&=\left.-2\pi\ii\frac{x-\ii}{(1-r)x-\ii(1+r)}\ee^{-\ii\lambda rx}\right|_{x=-\ii}\\ &=-2\pi\ii\ee^{-\lambda r}.
    \end{split}
\end{equation}
Hence 
\begin{equation}
    n(\lambda)=2\pi\int_0^1\,\dd r = 2\pi,\quad\lambda>0.
\end{equation}
We conclude that 
\begin{equation}
    \alpha(\lambda)=\frac{2\pi\ee^\lambda}{\mathrm{Ei}(2\lambda)+\ii\pi},\quad\lambda>0.
\end{equation}

This agrees precisely with the limit function $\widehat{\psi}(\lambda)$ that appears in Proposition~\ref{prop:weak-asymptotics} and Corollary~\ref{cor:weak-convergence}. To conclude the proof of Theorem~\ref{thm:minussoliton}, it is enough to 
establish the following lemma.

\begin{lemma}
Assume that $f\in\mathcal{C}(\R, L^2_+(\mathbb{R}))$ satisfies the following two conditions:
\begin{enumerate}
\item there is $\psi\in L^2_+(\mathbb{R})$ such that as $t\to+\infty$, 
\begin{equation}\label{hyp:weaklim}
\ee^{-\ii\pi/4}\sqrt{4\pi t}\ee^{\ii ty^2}f(t,2t \diamond)\rightharpoonup \widehat{\psi}(-\diamond)
\quad \text{in $L^2(-\infty,0)$};
\end{equation}
\item the $L^2$ norm of $f$ is preserved and equal to the $L^2$ norm of the weak limit: for every $t\in\R$,
\begin{equation}\label{hyp:norms}
\|f(t,\diamond)\|_{L^2(\mathbb{R})}=\|\psi\|_{L^2(\mathbb{R})}.
\end{equation}
\end{enumerate}
Then we have the following strong convergence in $L^2(\mathbb{R})$:
\begin{equation}
\lim_{t\to+\infty}\|f(t,\diamond)-\ee^{-\ii t\partial_{xx}}\psi\|_{L^2(\mathbb{R})}=0.
\end{equation} 
\end{lemma}
Here, $\ee^{-\ii t\partial_{xx}}$ denotes the Fourier multiplier $\ee^{\ii t\xi^2}$, and it defines the unitary flow of the linearized $\epsilon=0$ Benjamin-Ono equation $u_t =  \partial_x |D_x|u$ restricted to positive frequencies $\xi>0$, i.e., for initial data in $L^2_+(\mathbb{R})$.

Indeed, we apply the lemma to $f(t,\diamond)=\Pi u(t,\diamond)$. Condition~\eqref{hyp:weaklim} is guaranteed by Corollary~\ref{cor:weak-convergence} and condition~\eqref{hyp:norms} is a consequence of Corollary~\ref{spectralLuplus} because $\widehat{\psi}=\alpha=\widetilde{\Pi u_0}$. Writing $u=2\mathrm{Re}(\Pi u)$, we finally get Theorem~\ref{thm:minussoliton} with  
\begin{equation}
u_+=2\mathrm{Re}(\psi).
\end{equation}

\begin{proof}
Let $\psi^\sharp(t,x):=\ee^{\ii x^2/(4t)}\psi(x)$.  Then by expanding the norm squared,
\begin{multline}
    \|f(t,\diamond)-\ee^{-\ii t\partial_{xx}}\psi^\sharp(t,\diamond)\|^2_{L^2(\mathbb{R})} = 
    \|f(t,\diamond)\|^2_{L^2(\mathbb{R})} + \|\ee^{-\ii t\partial_{xx}}\psi^\sharp(t,\diamond)\|^2_{L^2(\mathbb{R})} \\{}-2\mathrm{Re}\left(\langle f(t,\diamond),\ee^{-\ii t\partial_{xx}}\psi^\sharp(t,\diamond)\rangle\right).
\end{multline}
Since $\ee^{-\ii t\partial_{xx}}$ and multiplication by $\ee^{\ii x^2/(4t)}$ are unitary, we apply the hypothesis \eqref{hyp:norms} to get
\begin{equation}
    \|f(t,\diamond)-\ee^{-\ii t\partial_{xx}}\psi^\sharp(t,\diamond)\|^2_{L^2(\mathbb{R})} =2\|\psi\|_{L^2(\mathbb{R})} -2\mathrm{Re}\left(\langle f(t,\diamond),\ee^{-\ii t\partial_{xx}}\psi^\sharp(t,\diamond)\rangle\right).
    \label{eq:L2diff-1}
\end{equation}
It is an exercise in Fourier theory to show that the following identity holds:
\begin{equation}
\ee^{-\ii t\partial_{xx}}\psi^\sharp(t,x) 
	=\frac{\ee^{\ii\pi/4}}{\sqrt{4\pi t}}\ee^{-\ii x^2/(4t)}\widehat{\psi}\left(-\frac{x}{2t}\right).
\end{equation}
Therefore,
\begin{equation}
\begin{split}
    \langle f(t,\diamond),\ee^{-\ii t\partial_{xx}}\psi^\sharp(t,\diamond)\rangle&=
    \int_\mathbb{R}\frac{\ee^{-\ii\pi/4}}{\sqrt{4\pi t}}f(t,x)\ee^{\ii x^2/(4t)}\widehat{\psi}\left(-\frac{x}{2t}\right)^*\,\dd x\\
    &= \frac{1}{2\pi}\int_\mathbb{R}\ee^{-\ii\pi/4}\sqrt{4\pi t}\ee^{\ii ty^2}f(t,2ty)\widehat{\psi}(-y)^*\,\dd y\\
    &=\frac{1}{2\pi}\langle \ee^{-\ii\pi/4}\sqrt{4\pi t}\ee^{\ii t\diamond^2}f(t,2t\diamond),\widehat{\psi}(-\diamond)\rangle.
\end{split}
\end{equation}
Given that $\psi\in L^2_+(\mathbb{R})$, the scalar product in the right-hand side is really the $L^2(-\infty,0)$ scalar product. Hence the hypothesis~\eqref{hyp:weaklim} is enough to conclude that
\begin{equation}
\lim_{t\to+\infty}\langle f(t,\diamond),\ee^{-\ii t\partial_{xx}}\psi^\sharp(t,\diamond)\rangle = 
    \frac{1}{2\pi}\langle\widehat{\psi},\widehat{\psi}\rangle=\langle\psi,\psi\rangle = \|\psi\|^2_{L^2(\mathbb{R})}.
\end{equation}
Using this in \eqref{eq:L2diff-1} gives
\begin{equation}
    \lim_{t\to+\infty}\|f(t,\diamond)-\ee^{-\ii t\partial_{xx}}\psi^\sharp(t,\diamond)\|_{L^2(\mathbb{R})} = 0.
\end{equation}
Then 
\begin{equation}
\begin{split}
    \|f(t,\diamond)-\ee^{-\ii t\partial_{xx}}\psi\|_{L^2(\mathbb{R})} &\le \|f(t,\diamond)-\ee^{-\ii t\partial_{xx}}\psi^\sharp(t,\diamond)\|_{L^2(\mathbb{R})} + \|\ee^{-\ii t\partial_{xx}}\psi-\ee^{-\ii t\partial_{xx}}\psi^\sharp(t,\diamond)\|_{L^2(\mathbb{R})}\\
    &=\|f(t,\diamond)-\ee^{-\ii t\partial_{xx}}\psi^\sharp(t,\diamond)\|_{L^2(\mathbb{R})} + \|\psi-\psi^\sharp(t,\diamond)\|_{L^2(\mathbb{R})}
\end{split}
\end{equation}
because $\ee^{-\ii t\partial_{xx}}$ is unitary. By dominated convergence, 
\begin{equation}
\lim_{t\to+\infty}\|\psi-\psi^\sharp(t,\diamond)\|_{L^2(\mathbb{R})}= 0,
\end{equation}
so the proof is complete.
\end{proof}

\appendix
    \section{The initial data \texorpdfstring{$u_0(x)=2/(x^2+1)$}{u0(x)=2/(x**2+1)}}\label{section:soliton}
To illustrate the use of Theorem~\ref{thm:inversion-formula}, we consider perhaps the simplest nontrivial case by choosing the initial condition to be
\begin{equation}
    u_0(x)=\frac{2}{x^2+1} = \frac{-\ii}{x-\ii} + \frac{\ii}{x+\ii}.
\label{eq:Lorentzian}
\end{equation}
Here $N=1$, $p_1=\ii$, and $c_1=-\ii$.  Hence $\ii c_1/\epsilon$ is a strictly positive integer $M$ whenever $\epsilon=1/M$.  Also, in this case it is known \cite{KodamaAblowitzSatsuma82,MillerW16a} that the reflection coefficient in the Fokas-Ablowitz inverse-scattering transform vanishes, and hence the corresponding solution of the Benjamin-Ono equation \eqref{eq:BO} is an ensemble of $M$ solitons (the $M$ Lax eigenvalues are simply related to the roots of the Laguerre polynomial of degree $M$, as was first shown in \cite{KodamaAblowitzSatsuma82}).  For $M=1$, the exact solution is known to be the Benjamin-Ono soliton:
\begin{equation}
u(t,x)=\frac{2}{(x-t)^2+1},\quad \epsilon=M^{-1}=1.
\label{eq:soliton}
\end{equation}
For larger integers $M$ a formula is known in the literature \cite{Matsuno79} for the corresponding solution. However that formula is in terms of $M\times M$ determinants.  We will now give a formula for the same solution as a $2\times 2$ determinant (because $N=1$ in \eqref{eq:lambda-formula}), that is moreover immediately reducible to a scalar expression.  Then we will set $M=1$ and show that the solution obtained is exactly \eqref{eq:soliton}.

\subsection{General \texorpdfstring{$M\in \mathbb{N}$}{}}
Letting $\epsilon=M^{-1}$ and taking $N=1$, $p_1=\ii$, and $c_1=-\ii$, the exponential factor appearing in the integrands of the matrix elements of $\mathbf{A}(t,x)$ and $\mathbf{B}(t,x)$ is single-valued because $M\in\mathbb{N}$:
\begin{equation}
    \ee^{-\ii h(z)/\epsilon}=\ee^{-\ii M(z-x)^2/(4t)}\frac{(z+\ii)^M}{(z-\ii)^M},
\end{equation}
according to \eqref{eq:h-def}.  Hence, from \eqref{eq:Bmatrix} and \eqref{eq:Btilde}, we have
\begin{equation}
    \det(\overline{\mathbf{B}}(t,x))=\ee^{\ii M(x-\ii)^2/(4t)}\det(\mathbf{B}(t,x))
\label{eq:detBtilde}
\end{equation}
where the matrix $\mathbf{B}(t,x)$ is
\begin{equation}
\mathbf{B}(t,x)=\begin{bmatrix}
    \displaystyle \int_{C_0}\ee^{-\ii M(z-x)^2/(4t)}\frac{(z+\ii)^M}{(z-\ii)^M}\,\dd z & \displaystyle \int_{C_0}\ee^{-\ii M(z-x)^2/(4t)}\frac{(z+\ii)^M}{(z-\ii)^{M+1}}\,\dd z\\
    \displaystyle \int_{C_1}\ee^{-\ii M(z-x)^2/(4t)}\frac{(z+\ii)^M}{(z-\ii)^M}\,\dd z & \displaystyle \int_{C_1}\ee^{-\ii M(z-x)^2/(4t)}\frac{(z+\ii)^M}{(z-\ii)^{M+1}}\,\dd z
\end{bmatrix}.
\label{eq:B-soliton}
\end{equation}
Note that since we are in the non-exceptional case for the contour $C_1$, it may be taken as a U-shaped contour surrounding only the pole at $z=\ii$ and tending to infinity in both directions with $\arg(z)=3\pi/4$.  Because the integrand is single-valued, it follows that the integrals over $C_1$ may be evaluated by residues.  It is convenient to first make the substitution $z=x+\ee^{-\ii\pi/4}w/\sigma$, with $\sigma:=\sqrt{M/(4t)}$ so that for $k=1,2$ 
\begin{equation}
\begin{split}
B_{2k}(t,x)&=\int_{C_1}\ee^{-\ii M(z-x)^2/(4t)}\frac{(z+\ii)^M}{(z-\ii)^{M+k-1}}\,\dd z \\ &= \left(\ee^{\ii\pi/4}\sigma\right)^{k-2}2\pi\ii\mathop{\mathrm{Res}}_{w=\ee^{-3\pi\ii/4}\sigma(x-\ii)}
\ee^{-w^2}\frac{\displaystyle(w-\ee^{-3\pi\ii/4}\sigma(x+\ii))^M}{(w-\ee^{-3\pi\ii/4}\sigma(x-\ii))^{M+k-1}} \\ &= \left.\frac{\displaystyle 2\pi\ii\left(\ee^{\ii\pi/4}\sigma\right)^{k-2}}{(M+k-2)!}\frac{\dd^{M+k-2}}{\dd w^{M+k-2}}\left[\ee^{-w^2}\left(w-\ee^{-3\pi\ii/4}\sigma(x+\ii)\right)^M\right]\right|_{w=\ee^{-3\pi\ii/4}\sigma(x-\ii)}.
\end{split}
\end{equation}
Expanding the derivative of the product, 
\begin{multline}
\frac{\dd^{M+k-2}}{\dd w^{M+k-2}}\left[\ee^{-w^2}\left(w-\ee^{-3\pi\ii/4}\sigma(x+\ii)\right)^M\right]\\
\begin{aligned}
    &=\sum_{n=0}^{M+k-2}
\frac{(M+k-2)!}{n!(M+k-2-n)!}\frac{\dd^n}{\dd w^n}\ee^{-w^2}\cdot\frac{\dd^{M+k-2-n}}{\dd w^{M+k-2-n}}\left(w-\ee^{-3\pi\ii/4}\sigma(x+\ii)\right)^M\\
&=\sum_{n=0}^{M+k-2}\frac{(M+k-2)!M!}{n!(M+k-2-n)!(n+2-k)!}\left(w-\ee^{-3\pi\ii/4}\sigma(x+\ii)\right)^{n+2-k}\frac{\dd^n}{\dd w^n}\ee^{-w^2}.
\end{aligned}
\end{multline}
The remaining derivative can be 
expressed in terms of Hermite polynomials $H_n(\diamond)$ by Rodrigues' formula \cite[Eqn.\@ 18.5.5]{DLMF}:
\begin{multline}
\frac{\dd^{M+k-2}}{\dd w^{M+k-2}}\left[\ee^{-w^2}\left(w-\ee^{-3\pi\ii/4}\sigma(x+\ii)\right)^M\right]\\
=\ee^{-w^2}\sum_{n=0}^{M+k-2}\frac{(M+k-2)!M!}{n!(M+k-2-n)!(n+2-k)!}\left(w-\ee^{-3\pi\ii/4}\sigma(x+\ii)\right)^{n+2-k}H_n(w).   
\end{multline}
Evaluating at $w=\ee^{-3\pi\ii/4}\sigma(x-\ii)$ then yields, for $k=1,2$,
\begin{equation}
    B_{2k}(t,x)=2\pi\ii M!(2\ii)^{2-k}\ee^{-\ii M(x-\ii)^2/(4t)}\sum_{n=0}^{M+k-2}
    \frac{\left(2\ii\ee^{-3\pi\ii/4}\sigma\right)^nH_n(\ee^{-3\pi\ii/4}\sigma(x-\ii))}{n!(M+k-2-n)!(n+2-k)!}.
    \label{eq:B-row2}
\end{equation}
Using \eqref{eq:detBtilde} and expanding $\det(\mathbf{B}(t,x))$ then cancels the exponential factor $\ee^{-\ii M(x-\ii)^2/(4t)}$ and leaves an expression for $\det(\overline{\mathbf{B}}(t,x))$ as a linear combination of two contour integrals:
\begin{equation}
    \det(\overline{\mathbf{B}}(t,x))=
    P_{2}(t,x)\int_{C_0}\ee^{-\ii (z-x)^2/(4t)}\frac{(z+\ii)^M}{(z-\ii)^M}\,\dd z - P_1(t,x)\int_{C_0}\ee^{-\ii (z-x)^2/(4t)}\frac{(z+\ii)^M}{(z-\ii)^{M+1}}\,\dd z,
\label{eq:det-Btilde-M}
\end{equation}
where $P_k(t,x):=\ee^{\ii M(x-\ii)^2/(4t)}B_{2k}(t,x)$ are polynomials in $x$ obtained from \eqref{eq:B-row2}.  According to \eqref{eq:tau-form}, the corresponding $M$-soliton solution of the Benjamin-Ono equation \eqref{eq:BO} is then computed as the logarithmic derivative of this scalar expression:
\begin{multline}
    u(t,x)=-2M^{-1}\mathrm{Im}\left(\frac{\partial}{\partial x}\log\left(
    P_{2}(t,x)\int_{C_0}\ee^{-\ii (z-x)^2/(4t)}\frac{(z+\ii)^M}{(z-\ii)^M}\,\dd z\right.\right.\\
    \left.\left.-P_1(t,x)\int_{C_0}\ee^{-\ii (z-x)^2/(4t)}\frac{(z+\ii)^M}{(z-\ii)^{M+1}}\,\dd z\right)\right).
\label{eq:u-N1-generalM}
\end{multline}
Here $C_0$ may be taken as any contour originating at $z=\infty\ee^{3\pi\ii/4}$ and terminating at $z=\infty\ee^{-\ii\pi/4}$ that passes below the pole $z=\ii$.  One can also replace $C_0$ with the real interval $[-R,R]$ and take the limit $R\to+\infty$ to define each integral as a convergent improper integral.
Equation \eqref{eq:u-N1-generalM} appears to be a previously-unknown explicit formula for the $M$-soliton solution with the initial data \eqref{eq:Lorentzian} and $\epsilon=M^{-1}$.

\subsection{The case \texorpdfstring{$M=1$}{}}
In the case $M=1$, using $H_0(w)=1$ and $H_1(w)=2w$, we easily calculate
\begin{equation}
    B_{21}(t,x)=-4\pi\ee^{-\ii (\ii-x)^2/(4t)},\quad B_{22}(t,x)=2\pi\ii\ee^{-\ii (\ii-x)^2/(4t)}\frac{t-x+\ii}{t},\quad M=1.
\end{equation}
It then follows from \eqref{eq:B-row2} and \eqref{eq:det-Btilde-M} for $M=1$ that
\begin{equation}
\begin{split}
    \det(\overline{\mathbf{B}}(t,x))&=\frac{2\pi\ii}{t}\int_{C_0}\ee^{-\ii (z-x)^2/(4t)}\left[ t-x+\ii +\frac{2\ii(\ii-x)}{z-\ii} +\frac{4t}{(z-\ii)^2}\right]\,\dd z\\
    &=\frac{2\pi\ii}{t}(t-(x-\ii))I_0(t,x) +\frac{4\pi}{t}(x-\ii)I_1(t,x)+8\pi\ii I_2(t,x),
    \end{split}
    \label{eq:detBtilde-N1}
\end{equation}
where, by the substitution $z=w+x$ (resulting in a translated contour $C_0'$),
\begin{equation}
    I_P(t,x):=\int_{C_0}\ee^{-\ii (z-x)^2/(4t)}\frac{\dd z}{(z-\ii)^P} = \int_{C_0'}\ee^{-\ii w^2/(4t)}\frac{\dd w}{(w+(x-\ii))^P},\quad P=0,1,2,\dots.
\end{equation}
A useful identity comes from integration by parts:
\begin{equation}
\begin{split}
    I_P(t,x)&=-\frac{1}{P-1}\int_{C_0'}\ee^{-\ii w^2/(4t)}\frac{\dd}{\dd w}\frac{1}{(w+(x-\ii))^{P-1}}\,\dd w \\
    &= -\frac{\ii}{2(P-1)t}\int_{C_0'}\ee^{-\ii w^2/(4t)}\frac{(w+(x-\ii))-(x-\ii)}{(w+(x-\ii))^{P-1}}\,\dd w \\
    &= \frac{\ii (x-\ii)}{2 (P-1)t}I_{P-1}(t,x)-\frac{\ii}{2(P-1)t}I_{P-2}(t,x),\quad P=2,3,4,\dots.
    \end{split}
    \label{eq:I-IBP}
\end{equation}
Using \eqref{eq:I-IBP} for $P=2$ allows \eqref{eq:detBtilde-N1} to be expressed in a simpler form as
\begin{equation}
    \det(\overline{\mathbf{B}}(t,x))=\frac{2\pi\ii}{t}(t-(x+\ii))I_0(t,x).
\end{equation}
Since the shifted contour $C_0'$ can be taken to be locally independent of $x$ by Cauchy's theorem, we see that
\begin{equation}
    \frac{\partial I_0}{\partial x}(t,x)=0.
\label{eq:I-deriv}
\end{equation}
Then, using \eqref{eq:I-deriv}, we get
\begin{equation}
    \frac{\partial}{\partial x}\det(\overline{\mathbf{B}}(t,x))=-\frac{2\pi\ii}{t}I_0(t,x).
\end{equation}
Consequently, according to \eqref{eq:tau-form}, the solution of the Benjamin-Ono equation in this case is (taking $\epsilon=M^{-1}$ for $M=1$)
\begin{equation}
    u(t,x)=-2\mathrm{Im}\left(\frac{\partial}{\partial x}\log(\det(\overline{\mathbf{B}}(t,x)))\right)=-2\mathrm{Im}\left(\frac{1}{x-t+\ii}\right) = \frac{2}{(x-t)^2+1}
\end{equation}
in agreement with \eqref{eq:soliton}.

\bibliographystyle{siamplain}
\bibliography{references}

\end{document}